\documentclass[MSNbibl,number,citesort,seceqn,dvips]{arxbj}
\usepackage{upgreek}
\usepackage{mathbh}
\usepackage{graphicx}

\aid{0}
\volume{20}
\issue{4}
\pubyear{2014}
\firstpage{1673}
\lastpage{1697}
\doi{10.3150/13-BEJ523} 

\makeatletter

\newcommand{\rmd}{\mathrm{d}}

\newcommand{\bigtimes}{\mathop{\mbox{\fontsize{17}{17}\selectfont{$\!\times$}}}}

\newcommand{\cal}{\mathcal}

\newcommand{\covM}{\bolds\Sigma}

\newcommand{\Var}{\operatorname{Var}}

\newcommand{\bbM}{\mathbb{M}}

\newcommand{\bo}{\mathbf o}
\renewcommand{\o}{\mathbf o}
\newcommand{\R}{{\mathbb R}}
\newcommand{\E}{{\mathbb E}}
\newcommand{\bbP}{\mathbb{P}}

\newcommand{\calA}{{\mathcal{A}}}
\newcommand{\calB}{{\cal B}}
\newcommand{\calH}{{\cal H}}

\newcommand{\ind}{\mathbh1}

\newcommand{\Lonk}{\,\displaystyle\mathop{\longrightarrow}^\mathrm{D}_{k \to\infty}\,}
\newcommand{\PasLonk}{\,\displaystyle\mathop{\longrightarrow}^{\mathbb{P}\mbox{-}\mathrm{a.s.}}_{k \to\infty}\,}
\newcommand{\Lond}{\stackrel{\mathrm{D}}{\longrightarrow}}
\newcommand{\eqd}{\stackrel{\mathrm{D}}{=}}

\newcommand{\lonr}{\,\displaystyle\mathop{\longrightarrow}_{r \to\infty}\,}

\newcommand{\timesop}{\mathop{\times}}

\newcommand{\X}{\mathbf{X}}

\newtheorem{Th}{Theorem}[section]
\newtheorem{Lemma}{Lemma}[section]
\newtheorem{Cor}[Th]{Corollary}
\newproclaim{Def}{Definition}[section]
\newproclaim{cond}{Condition}
\newremark{Remark}{Remark}

\makeatother

\begin{document}
\begin{frontmatter}

\title{Asymptotic goodness-of-fit tests for the Palm mark distribution of
stationary point processes with correlated marks}
\runtitle{Asymptotic goodness-of-fit tests for stationary point processes}

\begin{aug}
\author[1]{\inits{L.}\fnms{Lothar} \snm{Heinrich}\corref{}\thanksref{1}\ead[label=e1]{lothar.heinrich@math.uni-augsburg.de}},
\author[2]{\inits{S.}\fnms{Sebastian} \snm{L\"uck}\thanksref{2,e2}\ead[label=e2,mark]{sebastian.lueck@uni-ulm.de}} \and
\author[2]{\inits{V.}\fnms{Volker} \snm{Schmidt}\thanksref{2,e3}\ead[label=e3,mark]{volker.schmidt@uni-ulm.de}}
\address[1]{Institute of Mathematics, University of Augsburg, D-86135
Augsburg, Germany.\\ \printead{e1}}
\address[2]{Institute of Stochastics, Ulm University, D-89069 Ulm,
Germany.\\ \printead{e2};
\printead*{e3}}
\end{aug}

\received{\smonth{7} \syear{2012}}
\revised{\smonth{2} \syear{2013}}

%
\begin{abstract}
We consider spatially homogeneous marked point patterns in an
unboundedly expanding convex sampling window. Our main objective is to
identify the distribution of the typical mark by constructing an
asymptotic $\chi^2$-goodness-of-fit test. The corresponding test
statistic is based on a natural empirical version of the Palm mark
distribution and a smoothed covariance estimator which turns out to be
mean square consistent. Our approach does not require independent
marks and allows dependences between the mark field and the point
pattern. Instead we impose a suitable $\beta$-mixing condition on the
underlying stationary marked point process which can be checked for a
number of Poisson-based models and, in particular, in the case of
geostatistical marking. In order to study test performance, our test
approach is applied to detect anisotropy of specific Boolean models.
\end{abstract}

%

\end{frontmatter}

\section{Introduction}\label{secintint}

Marked point processes (MPPs) are versatile models for the statistical
analysis of data recorded at
irregularly scattered locations. The simplest marking scenario is
independent marking, where marks are given by a sequence of independent
and identically distributed random elements, which is also independent
of the underlying point pattern of locations.
A more complex class of models considers a so-called geostatistical
marking, where the marks are determined by the values of a random
field at the given locations. Although the random field usually
exhibits intrinsic spatial correlations, it
is assumed to be independent of the location point process (PP).
However, in many real datasets interactions between locations and marks occur.
Moreover, many marked point patterns arising in models from stochastic geometry
such as edge centers in (anisotropic) Voronoi-tessellations marked by
orientation or PPs marked by nearest-neighbour distances do not fit the
setting of geostatistical marking.
For recent asymptotic approaches to mark correlation analysis based on
mark variogram and mark covariance functions, we refer to \cite
{Guan04,Guan07b,HeMo12}.
The main goal of this paper is to investigate estimators of the Palm
mark distribution
$P_M^{\mathbf o}$ in point patterns exhibiting correlations between different
marks as well as between
marks and locations. The probability measure $P_M^{\mathbf o}$ can be
interpreted
as the distribution of the typical mark which denotes the mark of a
randomly chosen point of the pattern.
For any mark set $C$, we consider the scaled deviations $Z_k(C)=\sqrt{|W_k|}
((\widehat{P}{}^{\mathbf o}_M)_k(C) - P_M^{\mathbf
o}(C)
)$ as measure
of the distance between
$P_M^{\mathbf o}$ and an empirical Palm mark distribution $(\widehat
{P}{}^{\mathbf
o}_M)_k $.
In \cite{Hei12b}, we prove asymptotic normality of
the scaled deviation vector
$\mathbf{Z}_k = (Z_k(C_1),\ldots,Z_k(C_\ell))^T$
under appropriate strong mixing conditions
when the
observation window $W_k$ with volume $|W_k|$ grows unboundedly in all
directions as $k \to\infty$.
In this study, we in particular discuss consistent estimators for the
covariance matrix of the Gaussian limit of $\mathbf{Z}_k$. This
enables us to construct asymptotic $\chi^2$-goodness-of-fit tests for
the Palm mark distribution $P_M^{\mathbf o}$.
In a simulation study we apply our testing methodology to the
directional analysis
of random surfaces. 
For this purpose, we consider Cox processes on the boundary of Boolean models,
mark them with the local outer normal direction and test for a hypothetical
directional distribution. This allows to identify the rose of
directions of the surface
process associated with the Boolean model and represents an alternative
to a Monte Carlo test for the rose of direction suggested in \cite{benes01}.
The occurring MPPs differ fundamentally from the setting of
independent and geostatistical marking, for which functional central
limit theorems (CLTs) and
corresponding tests have been derived in \cite{Heinrich08,Pawlas09}.
In general, they also do not represent $m$-dependent MPPs.

Our paper is organized as follows. Section \ref{subbasnot}
introduces basic notation and
definitions. In Section \ref{secresults}, we present our main
results, which are proved
in Section \ref{secproofs}. In Section \ref{secexamples}, we
briefly discuss some
models satisfying the assumptions needed to prove our asymptotic results.
In the final Section \ref{secstaapp}, we study the performance of the
proposed tests by simulations.\vspace*{-3pt}

\section{Stationary marked point processes}\vspace*{-3pt}\label{subbasnot}

An MPP $\X_M= \sum_{n \ge1}\delta_{(X_n, M_n)}$ is a random locally
finite counting measure (see \cite{Daley0307}, Volume II, Chapter 9.1) on the
Borel sets of $\R^d \times\mathbb{M}$
with atoms $(X_n,M_n) $, where the mark space $\mathbb{M}$ is Polish
endowed with its Borel $\sigma$-algebra $\mathcal{B}(\mathbb{M})$.
Formally, $\X_M$ is a random element with values in the space $\mathsf
{N}_{\mathbb{M}}$ of locally finite counting measures $\varphi(\cdot
)$ on
${\mathcal B}(\R^d\times{\mathbb M})$, where $\mathsf{N}_{\mathbb
{M}}$ is equipped with the $\sigma$-algebra generated by all sets of
the form $\{\varphi\in{\mathsf N}_{{\mathbb M}}\dvt  \varphi(B \times
C)=j\}$
for $j \ge0$, bounded $B \in\mathcal{B}(\R^d) $, and $C \in
{\mathcal B}({\mathbb
M}) $.
Throughout we assume that $\X_M$ is simple, that is, all locations
$X_n$ in
$\R^d$ have multiplicity $1$ regardless which mark they have.
In what follows, we only consider \emph{stationary} MPPs,
which means that\looseness=-1
\[
\X_M \eqd\sum_{n\geq1}
\delta_{(X_n-x,M_n)} \qquad\mbox{for all } x \in\R^d.
\]\looseness=0
We always assume that the \emph{intensity} $\lambda=
\E\X_M([0,1)^d \times\mathbb{M})$ is finite.\vadjust{\goodbreak}

\subsection{Palm mark distribution}\label{secPalmMarkDistr}

For a stationary MPP $\X_M$ the probability measure
$P_M^{\o}$ on $\mathcal{B}(\mathbb{M})$ defined by
%
\begin{equation}
\label{defmardis} P_M^{\o}(C) = \frac{1}{\lambda}
\mathbb{E} \X_M\bigl([0,1)^d \times C\bigr),\qquad C\in
\mathcal{B}(\mathbb{M}),
\end{equation}
is called the \emph{Palm mark distribution} of $\X_M$. It can
be interpreted as the conditional distribution of the mark
of an atom of $\X_M$ located at the origin $\o$.
A random element $M_0$ in $\mathbb{M}$ with distribution $P_M^{\o}$ is
called \emph{typical mark} of $\X_M$.

\begin{Def}
An increasing sequence $\{W_k\}$ of convex and compact sets in $\R^d$
such that $\varrho(W_k)=\sup\{ r>0\dvt  B(x,r) \subset W_k$ for some $x
\in W_k\}
\rightarrow\infty$ as $k \to\infty$ is called a \emph{convex averaging
sequence} (briefly CAS). Here $ B(x,r)$ denotes the closed ball (w.r.t. the
Euclidean norm $\|\cdot\|$) with midpoint at $x\in\R^d$ and radius
$r \ge0 $.
\end{Def}
In the following, $|\cdot|$ denotes $d$-dimensional Lebesgue measure
and $\calH_{d-1}$ is the surface content
(i.e., $(d-1)$-dimensional Hausdorff measure).
Some results from convex geometry applied to CAS $\{W_k\}$ yield the following
inequalities (see \cite{Boehm04} and \cite{Heinrich08})
%
\begin{equation}
\label{ineqHeinrichPawlas} \frac{1}{\varrho(W_k)} \le\frac{\calH
_{d-1}(\partial W_k)}{|W_k|} \le
\frac{d}{\varrho(W_k)} \quad\mbox{and}\quad 1 - \frac{|W_k \cap(W_k-x)|}{|W_k|}
\le\frac{d \|x\|}{\varrho(W_k)}
\end{equation}
for $\|x\| \le\varrho(W_k) $.
Moreover, using the notation ${\overline H}_k = \{z \in\mathbb{Z}^d\dvt
|E_z \cap W_k|>0 \} $, where $E_z = [-1/2,1/2)^d +z$ for $z \in
\mathbb{Z}^d$, we have shown in
\cite{Hei12a,Hei12b} that for a CAS $\{W_k\}$
%
\begin{equation}\label{CASgridConvergence}
1\le\frac{\#{\overline H}_k}{|W_k|}\le1+ \frac
{|W_k\oplus B(o,\sqrt{d})|-|W_k|}{|W_k|}\mathop{\longrightarrow}
_{k \to\infty}1,
\end{equation}
which follows from Steiner's formula (see \cite{Schneider93}, page 197),
and (\ref{ineqHeinrichPawlas}).
If $\X_M$ is ergodic (for a precise definition see \cite{Daley0307},
Volume II, page 194), the individual ergodic theorem applied to MPPs (see
Theorem 12.2.IV and Corollary 12.2.V in
\cite{Daley0307}, Volume II) provides the $\mathbb{P}\mbox{-}\mathrm{a.s.}$ limits
%
\begin{equation}
\label{ergempmark} \widehat{\lambda}_k = \frac{\X_M(W_k \times
{\mathbb M})}{|W_k|}
\PasLonk\lambda\quad\mbox{and}\quad \bigl(\widehat{P}{}^{\o}_M
\bigr)_{k}(C) = \frac{\X_M(W_k \times
C)}{\X_M(W_k \times{\mathbb M})} \PasLonk
P^{\o}_M(C)
\end{equation}
for any $C \in\mathcal{B}(\mathbb{M})$ and an arbitrary CAS $\{W_k\}
$.

\subsection{Factorial moment measures and the covariance
measure}\label{secmomcummeas}

For any integer $m \ge1$, the $m$th \emph{factorial moment measure}
$\alpha_{\X_M}^{(m)}$ of
the MPP $\X_M$ is defined on $\mathcal{B}((\R^d\times{\mathbb
M})^m)$ by
%
\begin{equation}
\label{defalfisk} \alpha_{\X_M}^{(m)} \Biggl(
\bigtimes_{i=1}^m( B_i \times
C_i ) \Biggr) = \E{\sum^{\neq}_{n_1,\ldots,n_m
\geq1}}
\prod_{i=1}^m \bigl( \ind_{B_i}(X_{n_i})
\ind_{C_i}(M_{n_i}) \bigr),
\end{equation}
where the sum $\sum_{n_1,\ldots, n_m \geq1}^{\neq}$ runs over all
$m$-tuples of pairwise distinct indices $n_1,\ldots,n_m\ge1$ for
bounded $B_i \in\mathcal{B}(\R^d)$ and $C_i\in\mathcal{B}(\mathbb
{M}),
i=1,\ldots, m$. We also need
the $m$th factorial moment measure $\alpha_\X^{(m)}$ of the
\emph{unmarked} PP $\X(\cdot) = \X_M((\cdot)\times\bbM) = \sum_{n \ge
1}\delta_{X_n}(\cdot)$ defined on $\mathcal{B}((\R^d)^m)$ by
\[
\alpha_\X^{(m)} \Biggl( \bigtimes_{i=1}^mB_i
\Biggr) = \alpha_{\X_M}^{(m)} \Biggl( \bigtimes_{i=1}^m(
B_i \times\mathbb{M} ) \Biggr) \qquad\mbox{for bounded }
B_1,\ldots, B_m \in\mathcal{B}\bigl(\R^d\bigr).
\]
The stationarity of $\X_M$ implies that $\alpha_{\X}^{(m)}$ is
invariant under diagonal shifts,
which allows to define the $m$th \emph{reduced factorial moment measure}
$\alpha_{\X, \mathrm{red}}^{(m)}$
uniquely determined by the following disintegration formula
%
\begin{equation}
\label{defredalf} \alpha_\X^{(m)} \Biggl(
\bigtimes_{i=1}^mB_i \Biggr)=\lambda\int
_{B_1} \alpha_{\X,\mathrm{red}}^{(m)} \Biggl(
\bigtimes_{i=2}^m(B_i-x) \Biggr) \,\rmd x
\qquad\mbox{see \cite{Daley0307}, Volume II, Chapter 12.1}.
\end{equation}
The weak correlatedness between parts of $\X$ over distant Borel sets
may be expressed by the (factorial) \emph{covariance measure} 
$\gamma_\X^{(2)}$ on $\mathcal{B}((\R^d)^2)$ defined by
\[
\gamma_\X^{(2)} (B_1\times B_2 ) =
\alpha_\X^{(2)} (B_1 \times B_2 ) -
\lambda^2 |B_1| |B_2|.
\]
The \emph{reduced covariance measure} $\gamma^{(2)}_{\X,\mathrm{red}}\dvtx  \mathcal
{B}(\R^d)\to
[-\infty,\infty]$ is in general a
signed measure defined in analogy to (\ref{defredalf}) with $\gamma
_\X^{(2)}$ instead of $\alpha_\X^{(2)}$,
which shows that
\[
\gamma^{(2)}_{\X,\mathrm{red}}(B)=\alpha^{(2)}_{\X,\mathrm{red}}(B)
- \lambda|B| \qquad\mbox{for bounded } B\in\mathcal{B}\bigl(\R^d\bigr).
\]

\subsection{$m$-point Palm mark distribution}\label{subtwopoi}

For fixed mark sets $C_1,\ldots,C_m \in\mathcal{B}(\mathbb{M}),
m\ge1 $,
the $m$th factorial moment measure $\alpha_{\X_M}^{(m)}$ of the MPP (see
(\ref{defalfisk})) can be regarded as a measure on
$\mathcal{B}((\R^d)^m)$, which is
absolutely continuous w.r.t. $\alpha_\X^{(m)}$.
Thus, there exists a Radon--Nikodym density
$P_M^{x_1,\ldots, x_m}(C_1\times\cdots\times C_m)$,
such that for any $B_1,\ldots, B_m \in\mathcal{B}(\R^d)$,
%
\begin{equation}
\label{condRadonNikodym} \alpha_{\X_M}^{(m)} \Biggl(
\bigtimes_{i=1}^m(B_i\times C_i)
\Biggr) =\int_{\timesop_{i=1}^mB_i }P_M^{x_1,\ldots, x_m} \Biggl(
\bigtimes_{i=1}^mC_i \Biggr)
\alpha_\X^{(m)}\bigl(\rmd (x_1,\ldots,
x_m)\bigr).
\end{equation}
Since the mark space $\mathbb{M}$ is Polish, this Radon--Nikodym density
can be extended to a regular conditional distribution of the mark
vector $(M_1,\ldots, M_m)$ given that the corresponding atoms $X_1,\ldots, X_m$ are located at pairwise distinct points $x_1,\ldots,
x_m$, that is,
\[
P_M^{x_1,\ldots, x_m}( C ) = \mathbb{P}\bigl((M_1,\ldots,
M_m)\in C \mid X_1=x_1,\ldots,
X_m=x_m\bigr) \qquad\mbox{for } C \in\mathcal{B}\bigl(
\bbM^m\bigr).
\]
For details we refer to \cite{Kallenberg86}, page 164.
The above conditional distribution is called the \emph{$m$-point Palm
mark distribution} of $\X_M $.
In case of a stationary simple MPP $\X_M$, it is easily checked that
the one-point Palm mark distribution coincides
with the Palm mark distribution defined in (\ref{defmardis}).

The next result is indispensable to study asymptotic properties of
variance estimators for the empirical mark distribution. It extends a
formula stated in \cite{Hei10} for
unmarked PPs to the case of marked PPs. The proof of this extension
relies essentially on (\ref{condRadonNikodym}). Details are left to
the reader.

\begin{Lemma} \label{LemmaFormelMonster}
Let $\X_M = \sum_{n \ge1}\delta_{(X_n,M_n)}$ be an MPP satisfying
$\E\X_M(B \times{\mathbb M})^4 < \infty$ for all bounded $B \in
\mathcal{B}(\R^d)$, and let $f\dvtx  \R^d\times\R^d\times\mathbb
{M}^2 \mapsto\R^1$ be a Borel-measurable function such that the
second moment of $\sum_{p,q \geq1}^{\ne} | f(X_p, X_q,
M_p, M_q) |$ exists. Then,
%
\begin{eqnarray}
\label{streuung}
&&\Var\Biggl(\sum^{\neq}_{p,q \ge1}
f(X_p, X_q, M_p, M_q) \Biggr)
\nonumber\\
&&\quad= \int_{(\R^d)^2} \int_{{{\mathbb{M}}}^2}
f(x_1, x_2, u_1, u_2)
\bigl[f(x_1, x_2, u_1, u_2)+f(x_2,
x_1, u_2, u_1) \bigr] \nonumber\\
&&\qquad\hspace*{41pt}{}\times P_M^{x_1, x_2}
\bigl(\rmd (u_1,u_2) \bigr)\alpha_\X^{(2)}
\bigl(\rmd (x_1,x_2) \bigr)
\nonumber
\\
&&\qquad{}+ \int_{(\R^d)^3} \int_{{{\mathbb{M}}}^3}
f(x_1, x_2, u_1, u_2)
\bigl[f(x_1, x_3, u_1, u_3)+f(x_3,
x_1, u_3, u_1)
\nonumber
\\
&&\hspace*{125.1pt}\qquad{}+ f(x_2, x_3, u_2, u_3)+
f(x_3, x_2, u_3, u_2)
\bigr]\\
&&\qquad\hspace*{53.5pt}{}\times P_M^{x_1, x_2,
x_3} \bigl(\rmd (u_1,
u_2, u_3) \bigr)\alpha_\X^{(3)}
\bigl(\rmd (x_1,x_2, x_3) \bigr)
\nonumber
\\
&&\qquad{}+ \int_{(\R^d)^4}\int_{{{\mathbb{M}}}^4}
f(x_1, x_2, u_1, u_2)f(x_3,
x_4, u_3, u_4) \nonumber\\
&&\qquad\hspace*{54pt}{}\times\bigl[ P_M^{x_1, x_2, x_3, x_4}
\bigl(\rmd (u_1, u_2, u_3, u_4)
\bigr)\alpha_\X^{(4)} \bigl(\rmd (x_1,x_2,
x_3, x_4) \bigr)
\nonumber
\\
&&\hspace*{69pt}\qquad{}- P_M^{x_1, x_2} \bigl(\rmd (u_1,
u_2) \bigr)P_M^{x_3, x_4} \bigl(\rmd (u_3, u_4) \bigr) \alpha_\X^{(2)}
\bigl(\rmd (x_1,x_2) \bigr) \alpha_\X^{(2)}
\bigl(\rmd (x_3,x_4) \bigr) \bigr].
\nonumber
\end{eqnarray}
\end{Lemma}

\subsection{\texorpdfstring{$\beta$}{beta}-mixing coefficient and covariance inequality}\label{subbetmix}

For any $B \in\mathcal{B}(\R^d)$, let $\calA_{\X_M}(B)$ denote the
sub-$\sigma$-algebra of $\calA$ generated by the
restriction of the MPP $\X_M$ 
to the set $B \times{\mathbb M}$. For any $B, B^\prime\in\mathcal
{B}(\R^d)$, a natural
measure of dependence between $\calA_{\X_M}(B)$ and $\calA_{\X
_M}(B^\prime)$ can be
formulated in terms of the $\beta$-\emph{mixing} (or \emph{absolute
regularity}, respectively,
\emph{weak Bernoulli}) \emph{coefficient}
%
\begin{equation}
\label{defbetemm} \beta\bigl(\calA_{\X_M}(B),\calA_{\X_M}
\bigl(B^\prime\bigr) \bigr)= \frac{1}{2} 
\sup_{\{A_i\}, \{A'_j\}} \sum_{i,j}
\bigl| \bbP\bigl(A_i \cap A_j^{\prime}\bigr) -
\bbP(A_i) \bbP\bigl(A_j^{\prime}\bigr) \bigr|,
\end{equation}
where the supremum is taken over all finite partitions
$\{A_i\}$ and $\{A_j^{\prime}\}$ of $\Omega$ such that
$A_i \in\calA_{\X_M}(B)$ and $A_j^{\prime} \in\calA_{\X
_M}(B^\prime)$ for all $i,j $,
see \cite{Dou94} or \cite{Bra07} for a detailed discussion of this
and other mixing coefficients.
To quantify the degree of dependence of the MPP $\X_M$ on disjoint sets
$K_a = [-a,a]^d$ and
$K^c_{a+b}=\R^d \setminus K_{a+b}$, where $b \ge0$, we
introduce non-increasing rate functions
$\beta_{\X_M}^*, \beta_{\X_M}^{**}\dvtx [\frac{1}{2},\infty) \to
[0,\infty)$ depending on some constant $c_0\ge1$ such that
%
\begin{equation}
\label{betstaone} \beta\bigl(\calA_{\X_M}(K_a),
\calA_{\X_M}\bigl(K ^c_{a+b}\bigr) \bigr) \le
\cases{\beta_{\X_M}^*(b), &\quad for $\frac{1}{2} \le a \le
b/c_0$,
\vspace*{2pt}\cr
a^{d-1} \beta_{\X_M}^{**}(b),
&\quad for $\frac{1}{2} \le b/c_0 \le a$.}
\end{equation}
A stationary MPP $\X_M$ is called $\beta$-\emph{mixing} or \emph{absolutely
regular}, respectively, \emph{weak Bernoulli} if
both $\beta$-\emph{mixing rates} $\beta_{\X_M}^*(r)$ and $\beta_{\X
_M}^{**}(r)$ tend to $0$
as $r \to\infty$.
Note that any stationary $\beta$-mixing MPP $\X_M$ is mixing in the
usual sense and thus also
ergodic, see Lemma 12.3.II and Proposition~12.3.III in \cite
{Daley0307}, Volume II, page 206.
Our proofs of the asymptotic results in Section~\ref{secresults}
require at least polynomial decay of
$\beta_{\X_M}^*(r)$ and $\beta_{\X_M}^{**}(r)$ expressed by:
\renewcommand{\thecond}{$\bolds{\beta(\delta)}$}
\begin{cond}\label{condBeDe}
Let the MPP $\X_M$ satisfy (\ref{betstaone}) and
$\E\X_M( [0,1]^d \times\mathbb{M} )^{2+\delta} < \infty$ such that
\[
\int_1^\infty r^{d-1} \bigl(
\beta_{\X_M}^*(r) \bigr)^{\delta
/(2+\delta)} \,\rmd r < \infty\quad\mbox{and}\quad
r^{2d-1} \beta_{\X_M}^{**}(r) \lonr0 \qquad\mbox{for some }
\delta> 0.
\]
\end{cond}
A condition of this type based on (\ref{defbetemm}) and (\ref
{betstaone}) has been first verified for stationary \mbox{(Poisson-)}
Voronoi tessellations in \cite{Hei94}. It has proven adequate to
derive CLTs via Bernstein's blocking technique for spatial means
related with these tessellations observed in expanding cubic
observation windows. The proof of the below stated Theorem \ref{theasynor},
which is given in \cite{Hei12b}, extends Bernstein's method to
observation windows forming a CAS.
The following covariance bound in terms of the $\beta$-mixing coefficient
(\ref{defbetemm}) emerged first in \cite{Yosh76}, see also \cite{Bra07}.

\begin{Lemma}\label{lemcovbet} Let $Y$ and $Y^\prime$ denote the
restrictions of the MPP $\X_M$
to $B\times\mathbb{M} $ and $B^\prime\times\mathbb{M}$ for some
$B, B^\prime\in\mathcal{B}(\R^d) $, respectively.
Furthermore, let $\widetilde Y$ and $\widetilde Y^\prime$ be
independent copies of $Y$ and $Y^\prime$,
respectively. Then, for any ${\mathcal N}_{\mathbb{M}}\otimes
{\mathcal N}_{\mathbb{M}}$-measurable function
$f\dvtx \mathsf{N}_\mathbb{M}\times\mathsf{N}_\mathbb{M} \to[0,\infty
)$ and, for any $\eta> 0 $,
%
\begin{eqnarray}\label{covinebet}
&&
\bigl|\E f\bigl(Y,Y^\prime\bigr)-\E f\bigl(\widetilde Y,\widetilde
Y^\prime\bigr) \bigr|\nonumber\\
&&\quad\le 2 \beta\bigl({\calA_{\X_M}(B)},{
\calA_{\X_M}}\bigl(B^\prime\bigr)\bigr)^{{\eta
}/({1+\eta})}
\\
&&\qquad{}\times \max\bigl\{ \bigl(\E f^{1+\eta}\bigl(Y,Y^\prime\bigr)
\bigr)^{{1}/({1+\eta})}, \bigl(\E f^{1+\eta}\bigl(\widetilde
Y,\widetilde
Y^\prime\bigr) \bigr)^{{1}/({1+\eta})} \bigr\}.\nonumber
\end{eqnarray}
\end{Lemma}
If $f$ is bounded, then (\ref{covinebet}) remains valid for $\eta=
\infty$.

\section{Results}\label{secresults}


\subsection{Central limit theorem}

We consider a sequence of set-indexed empirical processes
$\{ Y_k(C), C \in\mathcal{B}(\mathbb{M}) \}$ defined by
%
\begin{eqnarray}
\label{ypskahtil} Y_k(C) &=& \frac{1}{\sqrt{|W_k|}}\sum
_{n \ge1} \ind_{W_k}(X_n) \bigl(
\ind_C(M_n) - P^{\o}_M(C)
\bigr)
\nonumber\\[-8pt]\\[-8pt]
&=& \sqrt{|W_k|}\, \widehat{\lambda}_k\, \bigl( \bigl(
\widehat{P}{}^{\o
}_M\bigr)_k(C)
- P_M^{\o}(C) \bigr),\nonumber
\end{eqnarray}
where $\{W_k\}$ is a CAS of observation windows in $\R^d$.
We will first state a multivariate CLT for the joint distribution of
$Y_k(C_1),\ldots,Y_k(C_\ell)$. For this,
let ``$\Lond$'' denote \emph{convergence in distribution} and
${\mathcal N}_{\ell}(a,\bolds{\Sigma})$
be an $\ell$-dimensional Gaussian vector with expectation (column)
vector $a \in\R^\ell$ and covariance matrix
$\bolds{\Sigma} = (\sigma_{ij})_{i,j=1}^\ell$.

\begin{Th}\label{theasynor}
Let $\X_M$ be a stationary MPP with $\lambda> 0$ satisfying Condition
\hyperref[condBeDe]{${\beta(\delta)}$}. Then
%
\begin{equation}
\label{asyequnor} \mathbf{Y}_k= \bigl(Y_k(C_1),\ldots,Y_k(C_\ell) \bigr)^\top\Lonk{\mathcal
N}_\ell(\o_{\ell},\bolds{\Sigma}) \qquad\mbox{for any }
C_1,\ldots,C_\ell\in\calB(\mathbb{M}),
\end{equation}
where $\o_{\ell}=(0,\ldots,0)^\top$ and the asymptotic covariance
matrix $\bolds{\Sigma} =
(\sigma_{ij})_{i,j=1}^\ell$ is given by the limits
%
\begin{equation}\label{covmattau}
\sigma_{ij}=\lim_{k \to\infty} \E Y_k(C_i)Y_k(C_j).
\end{equation}
\end{Th}

This CLT, which is proved in \cite{Hei12b} in detail, can be
reformulated for the empirical set-indexed process
$\{Z_k(C), C \in\mathcal{B}(\mathbb{M})\}$, where
\[
Z_k(C) = ( \widehat{\lambda}_k
)^{-1} Y_k(C)=\sqrt{|W_k|} \bigl( \bigl(
\widehat{P}{}^{\o}_M\bigr)_k(C)
- P^{\o}_M(C) \bigr).
\]
In other words, as refinement of the ergodic theorem (\ref{ergempmark}),
we derive asymptotic normality of a suitably scaled deviation of the
ratio-unbiased empirical Palm mark probabilities
$(\widehat{P}{}^{\o}_M)_k(C)$ from $P^{\o}_M(C)$
defined by
(\ref{defmardis}) for any $C \in\mathcal{B}(\mathbb{M}) $.
Since Condition \hyperref[condBeDe]{${\beta(\delta)}$} ensures the ergodicity of~$\X_M $,
the first limiting relation in (\ref{ergempmark}) combined with
Slutsky's lemma
yields the following result as a corollary of Theorem \ref{theasynor}.
%
\begin{Cor}
The conditions of Theorem \ref{theasynor} imply the CLT
\[
\mathbf{Z}_k = \bigl(Z_k(C_1),\ldots,Z_k(C_\ell)\bigr)^\top\Lonk{\mathcal
N}_\ell\bigl(\o_{\ell},\lambda^{-2} \bolds{\Sigma}
\bigr). \label{asyasynor}
\]
\end{Cor}

\subsection{\texorpdfstring{$\beta$}{beta}-mixing and integrability conditions}\label{secbetaintegrability}

In this subsection, we give a condition in terms of the mixing rate
$\beta_{\X_M}^*(r)$ which implies finite total variation of the
reduced covariance measure $\gamma^{(2)}_{\X,\mathrm{red}}$ and a certain
integrability condition (\ref{intconalf}) which expresses
weak dependence between any two marks located at far distant sites.
Both of these conditions enable us to show the unbiasedness,
respectively,
asymptotic unbiasedness of two estimators for the asymptotic
covariances (\ref{covmattau}).
Note that the total variation measure $|\gamma^{(2)}_{\X,\mathrm{red}}|$ of
$\gamma^{(2)}_{\X,\mathrm{red}}$
is defined as sum of the positive part $\gamma^{(2)+}_{\X,\mathrm{red}}$ and
negative part
$\gamma^{(2)-}_{\X,\mathrm{red}}$ of the Jordan decomposition of $\gamma
^{(2)}_{\X,\mathrm{red}}$,
that is,
\[
\gamma^{(2)}_{\X,\mathrm{red}}= \gamma^{(2)+}_{\X,\mathrm{red}}- \gamma^{(2)-}_{\X,\mathrm{red}}
\quad\mbox{and}\quad \bigl|\gamma^{(2)}_{\X,\mathrm{red}}\bigr| = \gamma^{(2)+}_{\X,\mathrm{red}} +
\gamma^{(2)-}_{\X,\mathrm{red}},
\]
where the positive measures $\gamma^{(2)+}_{\X,\mathrm{red}}$ and $\gamma
^{(2)-}_{\X,\mathrm{red}}$
are mutually singular, see \cite{Folland99}, page 87.

\begin{Lemma}\label{lemintgam}
Let $\X_M$ be a stationary MPP satisfying
\[
\E\X_M\bigl( [0,1]^d \times\mathbb{M}
\bigr)^{2+\delta} < \infty\quad\mbox{and}\quad \int_1^\infty
r^{d-1} \bigl(\beta_{\X_M}^*(r) \bigr)^{\delta
/(2+\delta)} \,\rmd r
< \infty\qquad\mbox{for some } \delta> 0
\]
with $\beta$-mixing rate $\beta_{\X_M}^*(r)$ defined in (\ref{betstaone}). 
Then
%
\begin{equation}
\label{intcummea} \bigl|\gamma_{\X, \mathrm{red}}^{(2)}\bigr|\bigl(
\R^d\bigr) < \infty
\end{equation}
and
%
\begin{equation}
\label{intconalf} \int_{\R^d} \bigl| P_M^{\o,x}(C_1
\times C_2) - P_M^{\o}(C_1)
P_M^{\o}(C_2) \bigr| \alpha^{(2)}_{\X,\mathrm{red}}(\rmd x) < \infty
\qquad\mbox{for any } C_1,C_2 \in\cal{B}(
\mathbb M).
\end{equation}
\end{Lemma}

\subsection{Representation of the asymptotic covariance matrix}
\label{seccovrep}

In Theorem \ref{theasynor}, we stated conditions for
asymptotic normality of the random vector
$\mathbf{Y}_k$. Clearly, (\ref{defmardis}) and (\ref{ypskahtil})
immediately imply that $\E Y_k(C) = 0$ for any $C\in\calB(\mathbb{M})$.
A representation formula for the asymptotic
covariance matrix ${\bolds\Sigma} $ is given in the following theorem.

\begin{Th}
\label{thereptau}
Let $\X_M$ be a stationary MPP satisfying (\ref{intconalf}) and let
$\{W_k\}$ be a CAS. Then the limits in (\ref{covmattau}) exist
and take the form
%
\begin{eqnarray}
\label{dartautil} \sigma_{ij} &=& \lambda\bigl(
P_M^{\o}(C_i \cap C_j) -
P_M^{\o}(C_i) P_M^{\o}(C_j)
\bigr) \nonumber\\
&&{}+ \lambda\int_{\R^d} \bigl( P_M^{\o,x}(C_i
\times C_j)
- P_M^{\o,x}(C_i \times\mathbb{M})
P_M^{\o}(C_j)\\
&&\hspace*{37pt}{} - P_M^{\o,x}(C_j
\times\mathbb{M}) P_M^{\o}(C_i) +
P_M^{\o}(C_i) P_M^{\o}(C_j)
\bigr) \alpha^{(2)}_{\X,\mathrm{red}}(\rmd x).
\nonumber
\end{eqnarray}

In particular, if $\X_M$ is marked independently, then
%
\begin{equation}
\label{tauindmar} \sigma_{ij}=\lambda\bigl( P_M^{\o}(
C_i \cap C_j )- P_M^{\o}(C_i)
P_M^{\o}(C_j) \bigr).
\end{equation}
\end{Th}

\subsection{Estimation of the asymptotic covariance matrix}\label{seccovest}

In Section \ref{secstaapp}, we will exploit the normal convergence
(\ref{asyequnor}) for
statistical inference of the typical mark distribution.
More precisely, assuming that the asymptotic covariance matrix
${\bolds\Sigma}$
is invertible, we consider asymptotic $\chi^2$-goodness-of-fit tests,
which are based on the distributional limit
%
\begin{equation}
\label{teststat} T_k = {\mathbf Y}_k^\top
\widehat{\covM}_k^{-1}{\mathbf Y}_k \Lonk\chi
^2_\ell,
\end{equation}
which is an immediate consequence of (\ref{asyasynor}) and Slutsky's
lemma, provided that
${\widehat{\bolds\Sigma}}_k$ is a consistent estimator for
${\bolds\Sigma}$.
As in (\ref{ypskahtil}), we use the notation ${\mathbf Y}_k= (
Y_k(C_1),\ldots, Y_k(C_\ell) )^\top$,
and the random variable $\chi^2_\ell$ is $\chi^2$-distributed with
$\ell$
degrees of freedom. In the following we will discuss several estimators for
${\bolds\Sigma}$. Our first observation is that the simple plug-in estimator
$\widehat{\bolds\Sigma}_k^{(0)}= (Y_k(C_i)Y_k(C_j)
)_{i,j=1}^{\ell}$ for ${\bolds\Sigma}$
is useless, since the determinant of $\widehat{\bolds\Sigma}_k^{(0)}$
vanishes. Instead of $\widehat{\bolds\Sigma}_k^{(0)}$ we take the
edge-corrected estimator $\widehat{\bolds\Sigma}_k^{(1)}=
((\widehat{\sigma}_{ij}^{(1)})_k )_{i,j=1}^{\ell}$ with
%
\begin{eqnarray}
\label{esttauone} \bigl(\widehat{\sigma}_{ij}^{(1)}
\bigr)_k&=& \frac{1}{|W_k|} \sum_{p \ge1}
\ind_{W_k}(X_p) \bigl( \ind_{C_i \cap
C_j}(M_p)-P_M^{\o}(C_i)
P_M^{\o}(C_j) \bigr)
\nonumber\\[-8pt]\\[-8pt]
&&{}+ \sum^{\neq}_{p,q \ge1} \frac{\ind_{W_k}(X_p)\ind_{W_k}(X_q) ( \ind
_{C_i}(M_p)-P_M^{\o}(C_i) )
(\ind_{C_j}(M_q)-P_M^{\o}(C_j) )}{|(W_k-X_p)\cap
(W_k-X_q)|}.
\nonumber
\end{eqnarray}
As an alternative, which can be implemented in a more efficient way, we
neglect the edge correction and consider the naive estimator $\widehat
{\bolds\Sigma}_k^{(2)}= ((\widehat{\sigma}_{ij}^{(2)})_k
)_{i,j=1}^{\ell}$ for
${\bolds\Sigma}$ with
\begin{eqnarray*}
\bigl(\widehat{\sigma}_{ij}^{(2)}\bigr)_k&=&
\frac{1}{|W_k|} \sum_{p \ge1} \ind_{W_k}(X_p)
\bigl( \ind_{C_i \cap
C_j}(M_p)-P_M^{\o}(C_i)
P_M^{\o}(C_j) \bigr)
\\
&&{}+ \frac{1}{|W_k|}\sum^{\neq}_{p,q\ge1}
\ind_{W_k}(X_p)\ind_{W_k}(X_q)
\bigl( \ind_{C_i}(M_p) - P_M^{\o}(C_i)
\bigr) \bigl( \ind_{C_j}(M_q) - P_M^{\o}(C_j)
\bigr).
\nonumber
\end{eqnarray*}
%
\begin{Th}\label{theconvar}
Let $\X_M$ be a stationary MPP satisfying (\ref{intconalf}) and let
$\{W_k\}$ be a CAS. Then $(\widehat{\sigma}_{ij}^{(1)})_k$ is an
unbiased estimator, whereas
$(\widehat{\sigma}_{ij}^{(2)})_k$
is an asymptotically unbiased estimator for $\sigma_{ij}$, where
$i,j=1,\ldots,\ell$.
\end{Th}
\begin{Remark*}
In general, neither $(\widehat{\sigma
}_{ij}^{(1)})_k$ nor $(\widehat{\sigma}_{ij}^{(2)})_k$ are
$L^2$-consistent estimators
for $\sigma_{ij}$, even if stronger moment and mixing conditions are imposed.
According to Lemma \ref{lemintgam}, the integrability condition
(\ref{intconalf}) in
Theorems \ref{thereptau} and \ref{theconvar} can be replaced by
the stronger Condition \hyperref[condBeDe]{${\beta(\delta)}$}.
In order to obtain an $L^2$-consistent estimator,
we introduce a smoothed version of the unbiased estimator in (\ref
{esttauone}), which is based on some kernel function and a sequence
of bandwidths depending on the CAS $\{W_k\}$.
\end{Remark*}
\renewcommand{\thecond}{$\bolds{(wb)}$}
\begin{cond}\label{condWB}
Let $w\dvtx  \R\mapsto\R$ be a non-negative,
symmetric,
Borel-measurable \emph{kernel function} satisfying $w(x)\longrightarrow
w(0)=1$ as $x \to0 $. In addition, assume that $w(\cdot)$ is bounded by
$m_w<\infty$ and vanishes outside $B({\mathbf o},r_w) $ for some $r_w\in
(0,\infty)$.
Further, associated with $w(\cdot)$ and some given CAS $\{W_k\}$, let
$\{b_k\}$ be a
sequence of positive \emph{bandwidths} such that
%
\begin{equation}
\label{condbk1} \frac{\varrho(W_k)}{2 d r_w |W_k|^{1/d}}\geq b_k
\mathop{
\longrightarrow} _{k \to\infty}0,\qquad %
b_k^d
|W_k| \mathop{\longrightarrow} _{k \to\infty}\infty\quad\mbox{and}\quad
b_k^{{3}d/{2}} |W_k| \mathop{\longrightarrow}
_{k \to
\infty}0.
\end{equation}
\end{cond}
%
\begin{Th}
\label{thmestTauThree}
Let $\{W_k\}$ be an arbitrary CAS and $w(\cdot)$ be a
kernel function with an associated sequence of bandwidths $\{b_k\}$ satisfying
Condition \hyperref[condWB]{$(wb)$}. If the stationary MPP $\X_M$ satisfies
%
\begin{equation}\label{condConsistencyTau}
\E\X_M\bigl( [0,1]^d \times\mathbb{M}
\bigr)^{4+\delta} < \infty\quad\mbox{and}\quad \int_1^\infty
r^{d-1} \bigl(\beta_{\X_M}^*(r) \bigr)^{\delta
/(4+\delta)} \,\rmd r
< \infty
\end{equation}
for some $\delta>0$ with $\beta$-mixing rate $\beta_{\X_M}^*(r)$
defined in (\ref{betstaone}), then 
$\E( \sigma_{ij} - (\widehat{\sigma}_{ij}^{(3)})_k )^2
\mathop{\longrightarrow
}\limits_{k \to\infty}0 $,
where $(\widehat{\sigma}_{ij}^{(3)})_k$ is a smoothed covariance
estimator defined by
\begin{eqnarray*}
\bigl(\widehat{\sigma}_{ij}^{(3)}
\bigr)_k &=& \frac{1}{|W_k|} \sum_{p \ge1}
\ind_{W_k}(X_p) \bigl( \ind_{C_i \cap C_j}(M_p)
- P^{\o}_M(C_i) P^{\o}_M(C_j)
\bigr)
\\
&&{} + \sum^{\neq}_{p,q \ge1}
\frac{\ind_{W_k}(X_p) \ind_{W_k}(X_q) ( \ind_{C_i}(M_p) -
P_M^{\o}(C_i) )
( \ind_{C_j}(M_q) - P_M^{\o}(C_j) )}{|(W_k-X_p) \cap
(W_k-X_q)|} \\
&&\qquad\hspace*{12.5pt}{}\times w \biggl(\frac{\|X_q-X_p\|}{b_k |W_k|^{1/d}} \biggr).
\nonumber
\end{eqnarray*}
\end{Th}
\begin{Remark*}
The full strength of condition (\ref{condConsistencyTau}) imposed on the
$\beta$-mixing rate $\beta_{\X_M}^*(r)$ introduced in (\ref{betstaone})
is only needed to prove the consistency result of
Theorem \ref{thmestTauThree}.
In order to prove (\ref{intcummea}), (\ref{intconalf}), and
Theorem \ref{theasynor} it suffices to take the
somewhat smaller non-increasing rate function
%
\begin{equation}\label{condbetaDeltaA}
\beta_{\X_M}^*(r) = \beta\bigl(\calA_{\X_M}(K_a),
\calA_{\X
_M}\bigl(K_{a+r}^c\bigr) \bigr) \qquad\mbox{for }
r \ge a = 1/2.
\end{equation}

Moreover, as shown in \cite{Hei12a}, the assertions of Theorems \ref
{theasynor} and \ref{thereptau} remain valid if in
Condition \hyperref[condBeDe]{${\beta(\delta)}$} the rate functions $\beta_{\X_M}^*$ and
$\beta_{\X_M}^{**}$ (defined by the $\beta$-mixing coefficient
(\ref{defbetemm})) are replaced by the corresponding rate functions
derived as in (\ref{betstaone}) from the smaller
$\alpha$-mixing coefficient
\[
\alpha\bigl(\calA_{\X_M}(B),\calA_{\X_M}
\bigl(B^\prime\bigr) \bigr)= \sup\bigl\{ \bigl|\bbP\bigl(A\cap
A^\prime\bigr)-\bbP(A)\bbP\bigl(A^\prime\bigr) \bigr|\dvt  A\in
\calA_{\X_{M}}(B),A^\prime\in\calA_{\X_{M}}
\bigl(B^\prime\bigr)\bigr\},
\]
which results in a slightly weaker mixing condition on $\X_M$, see
\cite{Bra07} for a comparison of $\alpha$- and $\beta$-mixing. A
covariance inequality for the $\alpha$-mixing case similar to
(\ref{covinebet}) can be found in~\cite{Dou94}, see \cite{Hei12a} for
an improved version. Since for most of the MPP models the subtle
differences between $\alpha$- and $\beta$-mixing are irrelevant we
present our results under the unified assumptions of
Condition~\hyperref[condBeDe]{${\beta(\delta)}$} and
(\ref{condConsistencyTau}) with $\beta$-mixing rate functions as
defined in (\ref{betstaone}).

Concerning the shape of the observation windows $\{W_k\}$, the
relations (\ref{ineqHeinrichPawlas}) and (\ref{CASgridConvergence})
are essential in the proofs of our results. However, there exist
sequences of not necessarily convex sets $\{W_k\}$ which satisfy
(\ref{ineqHeinrichPawlas}) and (\ref{CASgridConvergence}), see
references in \cite{Hei12a}.
\end{Remark*}

\section{Proofs}
\label{secproofs}

\subsection{Proof of Lemma \texorpdfstring{\protect\ref{lemintgam}}{3.1}}

By definition of the signed measures $\gamma^{(2)}_\X$ and
$\gamma^{(2)}_{\X,\mathrm{red}}$ in Section \ref{secmomcummeas} and
using algebraic
induction, for any bounded Borel-measurable function $g\dvtx  (\R^d)^2 \to
\R^1 $ we obtain the relation
%
\begin{equation}
\label{defgamunm} \lambda\int_{\R^d}\int
_{\R^d} g(x,y) \gamma^{(2)}_{\X,\mathrm{red}}(\rmd y)
\,\rmd x = \int_{(\R^d)^2} g(x,y-x) \gamma^{(2)}_\X
\bigl(\rmd (x,y)\bigr).
\end{equation}
Let $H^+, H^-$ be a Hahn decomposition of $\R^d$ for $\gamma
^{(2)}_{\X,\mathrm{red}}$, that is,
\[
\gamma^{(2)+}_{\X,\mathrm{red}}(\cdot)= \gamma^{(2)}_{\X,\mathrm{red}}
\bigl(H^+ \cap(\cdot)\bigr) \quad\mbox{and}\quad \gamma^{(2)-}_{\X,\mathrm{red}}(
\cdot)= -\gamma^{(2)}_{\X,\mathrm{red}}\bigl(H^-\cap(\cdot)\bigr).
\]
We now apply (\ref{defgamunm}) for $g(x,y) = \ind_{E_{\o}}(x)
\ind_{H^+\cap E_z}(y) $,
where
$E_z =
[-\frac{1}{2},\frac{1}{2})^d + z$ for $z\in\mathbb{Z}^d $.
Combining this with
the definition (\ref{defredalf}) of the (reduced) second factorial
moment measures
$\alpha^{(2)}_\X$ and $\alpha^{(2)}_{\X,\mathrm{red}}$ of the unmarked PP
$\X= \sum_{i\ge1} \delta_{X_i} $
and using the relation
\[
\gamma_\X^{(2)}(A \times B) = \alpha_\X^{(2)}(A
\times B) - \lambda^2 |A| |B| \qquad\mbox{for all bounded } A,B\in
\mathcal{B}\bigl(\R^d\bigr),
\]
we obtain
\begin{eqnarray*}
\lambda\gamma^{(2)}_{\X,\mathrm{red}}\bigl(H^+\cap E_z
\bigr) &=& \int_{(\R^{d})^2} \ind_{E_{\o}}(x)
\ind_{H^+ \cap
E_z}(y-x)\alpha^{(2)}_\X\bigl(\rmd  (x,y)
\bigr) - \lambda^2 |E_{\o}| \bigl|H^+ \cap E_z\bigr|
\\
&=& \E\sum^{\neq}_{i,j\ge1}
\ind_{E_{\o}}(X_i)\ind_{H^+\cap
E_z}(X_j-X_i)
- \E\X(E_{\o}) \E\X\bigl(H^ + \cap E_z\bigr).
\end{eqnarray*}
Since $\o\notin H^+ \cap E_z$ for $z\in\mathbb{Z}^d$ with $|z|\ge2$ we
may continue with
%
\begin{eqnarray}\label{gamequ1}
\lambda\gamma^{(2)}_{\X,\mathrm{red}}\bigl(H^+ \cap E_z
\bigr) &=& \E\sum_{i\ge1}\delta_{X_i}(E_{\o})
\X\bigl(\bigl(H^+ \cap E_z\bigr) + X_i \bigr) - \E
\X(E_{\o}) \E\X\bigl(H^+\cap E_z\bigr)
\nonumber\\[-8pt]\\[-8pt]
&=& \E f\bigl(Y,Y^\prime_z\bigr) - \E f\bigl(\widetilde
Y,\widetilde Y^\prime_z\bigr) \qquad\mbox{for } |z| \ge2,\nonumber
\end{eqnarray}
where
%
\begin{equation}
\label{gamequ2} f\bigl(Y,Y^\prime_z\bigr) = \sum
_{i\ge1}\delta_{X_i}(E_{\o}) \X
\bigl(\bigl(H^+ \cap E_z\bigr) + X_i \bigr) \le\X
(E_{\o} ) \X(E_z\oplus E_{\o} )
\end{equation}
with $Y(\cdot) = \sum_{i \ge1}\delta_{X_i} ((\cdot)\cap E_{\o
} )$, respectively,
$Y^\prime_z(\cdot)=\sum_{j\ge1}\delta_{X_j} ((\cdot) \cap
(E_z\oplus E_{\o}) )$
being restrictions of the stationary PP $\X= \sum_{i \ge1}\delta_{X_i}$
to $E_{\o}$, respectively, $E_z\oplus E_{\o} = [-1,1)^d+z $. Further, let
$\widetilde Y$
and $\widetilde Y^\prime_z$ denote copies of the PPs $Y$ and
$Y^\prime_z $, respectively, which are assumed to be independent
implying that
$\E f(\widetilde Y, \widetilde Y'_z) = \E\X(E_{\o}) \E\X(H^ +
\cap E_z) $.
Since $Y$ is measurable w.r.t. $\mathcal{A}_\X(E_{\o})$, whereas
$Y^\prime_z$ is
$\mathcal{A}_\X(\R^d\setminus[-(|z|-1),|z|-1]^d )$-measurable,
we are in a position to apply Lemma~\ref{lemcovbet} with
$\beta(\mathcal{A}_\X(E_{\o}),\mathcal{A}_\X(\R^d
\setminus[-(|z|-1),|z|-1]^d )
\le\beta_{\X_M}^*(|z|-\frac{3}{2})$\vspace*{1pt} for $|z|\ge(c_0+3)/2 \ge2 $.
Hence, the estimate
(\ref{covinebet}) together with (\ref{gamequ1}) and (\ref
{gamequ2}) yields
\[
\bigl| \lambda\gamma^{(2)}_{\X,\mathrm{red}}\bigl(H^+\cap E_z
\bigr) \bigr| \le2 \biggl( \beta^*_{\X_M}\biggl(|z|-\frac{3}{2}\biggr)
\biggr)^{{\eta
}/({1+\eta})} \bigl( \max\bigl\{ \E f^{1+\eta}
\bigl(Y,Y^\prime_z\bigr), \E f^{1+\eta}\bigl(
\widetilde Y,\widetilde Y^\prime_z\bigr) \bigr\}
\bigr)^{{1}/({1+\eta})},
\]
where the maximum term on the rhs has the finite upper bound
$2^{d(1+\eta)} \E\X(E_{\o})^{2+2\eta}$ for $\delta= 2 \eta> 0$
in accordance with our assumptions. This is seen from (\ref{gamequ2})
using the Cauchy--Schwarz inequality and the stationarity of $\X$
giving
\[
\E f^{1+\eta}\bigl(Y, Y'_z\bigr) \le\bigl( \E
\X(E_{\o})^{2+2\eta} \E\X\bigl([-1,1]^d
\bigr)^{2+2\eta} \bigr)^{1/2} \le2^{d(1+\eta)} \E\X
(E_{\o})^{2+2\eta}
\]
and the same upper bound for $\E f^{1+\eta}(\widetilde Y,\widetilde
Y^\prime_z) $. By combining all the above estimates with
$\lambda\gamma^{(2)}_{\X,\mathrm{red}}(H^+ \cap[-\frac{3}{2},\frac{3}{2})^d)
\le3^d \E\X(E_{\o})^2$, we arrive at
\begin{eqnarray*}
&&
\lambda\gamma^{(2)}_{\X,\mathrm{red}}\bigl(H^+\bigr) \\
&&\quad\le3^d
\E\X(E_{\o})^2 + 2^{d+1} \bigl(\E
\X(E_{\o})^{2+\delta} \bigr)^{{2}/({2+\delta})} \sum
_{z\in\mathbb{Z}^d:  |z| \ge(c_0+3)/2} \biggl( \beta^*_{\X_M}\biggl(|z|-
\frac{3}{2}\biggr) \biggr)^{{\delta
}/({2+\delta})}.
\end{eqnarray*}
By the assumptions of Lemma \ref{lemintgam} the moments and the
series on the rhs
are finite and the same bound can be derived for
$-\lambda\gamma^{(2)}_{\X,\mathrm{red}}(H^-)$ which shows the validity of
(\ref{intcummea}).

The proof of (\ref{intconalf}) resembles that of (\ref
{intcummea}). First, we extend
the identity (\ref{defgamunm}) to the (reduced) second factorial
moment measure of the
MPP $\X_M$ defined by (\ref{defalfisk}) and (\ref
{condRadonNikodym}) for $m=2$
which reads as follows:
\begin{eqnarray}
&&\lambda\int_{\R^d} \int_{\R^d} g(x,y)
P_M^{\o,x}(C_1 \times C_2)
\alpha_{\X, \mathrm{red}}^{(2)}(\rmd y)\,\rmd x
\nonumber
\\
&&\quad= \int_{(\R^d)^2} g(x,y-x) P_M^{x,y}(C_1
\times C_2) \alpha_\X^{(2)} \bigl(\rmd (x,y)
\bigr)
\nonumber
\\
&&\quad= \E\sum^{\neq}_{i,j \ge1}
g(X_i,X_j-X_i) \ind_{C_1}(M_i)
\ind_{C_2}(M_j).
\nonumber
\end{eqnarray}
For the disjoint Borel sets $G^+$ and $G^-$ defined by
\[
G^{+(-)} = \bigl\{ x \in\R^d\dvt  P_M^{\o,x}(C_1
\times C_2) \ge(<) P_M^{\o}(C_1)
P_M^{\o}(C_2) \bigr\}
\]
we replace $g(x,y)$ in the above relation by $g^{\pm}(x,y) =
\ind_{E_{\o}}(x) \ind_{E_z^{\pm}}(y) $, where $E_z^{\pm} =
G^{\pm}
\cap E_z$ for $|z| \ge2 $, and consider the
restricted MPPs
$Y_{\o}(\cdot) = \X_M ((\cdot)\cap(E_{\o}\times C_1) )$,
$Y'_{z,\pm}(\cdot) = \X_M ((\cdot)\cap((E_z^{\pm} \oplus
E_{\o})\times C_2) )$ and their
copies $\widetilde Y_{\o}$ and $\widetilde Y'_{z,\pm} $, which are
assumed to be stochastically
independent. Further, in analogy to (\ref{gamequ2}), define
\[
f\bigl(Y_{\o},Y'_{z,\pm}\bigr) = \sum
_{i \ge1}\delta_{(X_i,M_i)}(E_{\o
}\times
C_1) \X_M \bigl( \bigl(E_z^{\pm}
+ X_i\bigr) \times C_2 \bigr) \le\X(E_{\o} )
\X( E_z \oplus E_{\o} ).
\]

It is rapidly seen that for $|z| \ge2 $
\[
\E f\bigl(Y_{\o},Y'_{z,\pm}\bigr) = \lambda
\int_{E_z^{\pm}} P_M^{\o,x}(C_1
\times C_2) \alpha^{(2)}_{\X,\mathrm{red}}(\rmd x)
\]
and
\[
\E f\bigl(\widetilde Y_{\o},\widetilde Y'_{z,\pm}
\bigr) = \E\X_M(E_{\o}\times C_1) \E
\X_M\bigl(E_z^{\pm} \times C_2
\bigr) = \lambda^2 P_M^{\o}(C_1)
P_M^{\o}(C_2) \bigl|E_z^{\pm}\bigr|
\]
and in the same way as in the foregoing proof we find that, for $|z|\ge
(c_0+3)/2 $,
\[
\bigl| \E f\bigl(Y_{\o},Y'_{z,\pm}\bigr) - \E f\bigl(
\widetilde Y_{\o},\widetilde Y'_{z,\pm}\bigr) \bigr|
\le2^{d+1} \bigl(\E\X(E_{\o})^{2+\delta}
\bigr)^{
{2}/({2+\delta})} \biggl( \beta^*_{\X_M}\biggl(|z|-\frac{3}{2}
\biggr) \biggr)^{{\delta
}/({2+\delta})}.
\]
Finally, the decomposition $\alpha_{\X, \mathrm{red}}^{(2)}(\cdot) = \gamma
_{\X,
\mathrm{red}}^{(2)}(\cdot) + \lambda|\cdot|$ together with the previous estimate
leads to
\begin{eqnarray*}
&&
\lambda\int_{E_z} \bigl| P_M^{\o,x}(C_1
\times C_2) - P_M^{\o}(C_1)
P_M^{\o}(C_2) \bigr| \alpha^{(2)}_{\X,\mathrm{red}}(
\rmd x) \\
&&\quad= \E f\bigl(Y_{\mathbf o},Y'_{z,+}\bigr) - \E
f\bigl(\widetilde Y_{\mathbf o},\widetilde Y'_{z,+}
\bigr)
\\
&&\qquad{} - \bigl( \E f\bigl(Y_{\mathbf o},Y'_{z,-}\bigr) -
\E f\bigl(\widetilde Y_{\mathbf o},\widetilde Y'_{z,-}
\bigr) \bigr) - \lambda P_M^{\o}(C_1)
P_M^{\o}(C_2) \bigl(\gamma^{(2)}_{\X,\mathrm{red}}
\bigl( E_z^+ \bigr) - \gamma^{(2)}_{\X,\mathrm{red}}\bigl(
E_z^- \bigr) \bigr)
\\
&&\quad\le 2^{d+2} \bigl(\E\X(E_{\o})^{2+\delta}
\bigr)^{
{2}/({2+\delta})} \biggl(\beta_{\X_M}^*\biggl(|z|-\frac{3}{2}
\biggr) \biggr)^{{\delta
}/({2+\delta})} \\
&&\qquad{}+ \lambda\bigl|\gamma^{(2)}_{\X, \mathrm{red}}\bigr|(E_z)
\qquad\mbox{for } |z| \ge(c_0+3)/2.
\end{eqnarray*}
Thus, the sum over all $z \in\mathbb{Z}^d$
is finite in view of our assumptions and the above-proved relation
(\ref{intcummea}) which completes the proof of Lemma \ref{lemintgam}.

\subsection{Proof of Theorem \texorpdfstring{\protect\ref{thereptau}}{3.3}}

It suffices to show (\ref{dartautil}), since independent marks
imply that $P_M^{{\mathbf o},x}(C_1\times C_2)
= P^{\mathbf o}_M(C_1) P^{\mathbf o}_M(C_2)$ for $x \ne{\mathbf o}$ and any
$C_1,C_2 \in\calB(\mathbb M)$ so that the integrand on the rhs of
(\ref{dartautil}) disappears which yields
(\ref{tauindmar}) for stationary independently MPPs.
By the very definition of $Y_k(C)$, we obtain that
%
\begin{eqnarray}\label{streu2}
&&
\operatorname{Cov} \bigl( Y_k(C_i), Y_k(C_j)
\bigr) \nonumber\\
&&\quad= \frac{1}{|W_k|} \E\sum_{p\ge1}
\ind_{W_k}(X_p) \bigl(\ind_{C_i}(M_p)-P_M^{\mathbf o}(C_i)
\bigr) \bigl(\ind_{C_j}(M_p)-P_M^{\mathbf o}(C_j)
\bigr)
\\
&&\qquad{}+ \frac{1}{|W_k|} \E\sum_{p,q \geq1}^{\neq}
\ind_{W_k}(X_p)\ind_{W_k}(X_q)
\bigl(\ind_{C_i}(M_p)-P_M^{\mathbf
o}(C_i)
\bigr) \bigl(\ind_{C_j}(M_q)-P_M^{\mathbf o}(C_j)
\bigr).\nonumber
\end{eqnarray}
Expanding the difference terms in the parentheses leads to eight
expressions which, up to constant factors, take either the form
\[
\E\sum_{p \ge1} \ind_{W_k}(X_p)
\ind_C(M_p) = \lambda|W_k|
P^{\bo}_M(C)\vadjust{\goodbreak}
\]
or
\begin{eqnarray*}
&&\E\sum
_{p,q \geq
1}^{\neq} \ind_{W_k}(X_p)
\ind_{W_k}(X_q)\ind_{C_i}(M_p)
\ind_{C_j}(M_q)
\\
&&\quad= \int_{(\R^d)^2}\ind_{W_k}(x)\ind_{W_k}(y)
P_M^{\bo,y-x}(C_i \times C_j)
\alpha_\X^{(2)}\bigl(\rmd (x,y)\bigr) \\
&&\quad= \lambda\int
_{{\R}^d}P_M^{\bo,y}(C_i
\times C_j) \gamma_k(y) \alpha_{\X,\mathrm{red}}^{(2)}(\rmd y),
\end{eqnarray*}
where $y \mapsto\gamma_k(y)=|W_k \cap(W_k-y)|$ denotes the set covariance
function of $W_k $. Summarizing all these terms gives
\begin{eqnarray*}
&&
\operatorname{Cov} \bigl(Y_k(C_i),Y_k(C_j)
\bigr) \\
&&\quad= \lambda\bigl(P_M^{\mathbf o}(C_i \cap
C_j)-P_M^{\mathbf
o}(C_i)P_M^{\mathbf o}(C_j)
\bigr) \\
&&\qquad{}+ \lambda\int_{{\mathbb R}^d} \frac{\gamma_k(x)}{|W_k|} \bigl(
P_M^{\bo,x}(C_i \times C_j)
- P_M^{\mathbf o}(C_i) P_M^{\bo,x}(C_j
\times\mathbb{M})\\
&&\hspace*{62.5pt}\qquad{} -P_M^{\mathbf o}(C_j)
P_M^{\bo,x}(C_i \times\mathbb{M}) +
P_M^{\mathbf o}(C_i) P_M^{\mathbf o}(C_j)
\bigr) \alpha_{\X,\mathrm{red}}^{(2)}(\rmd x).
\end{eqnarray*}
The integrand in the latter formula is dominated by the sum
\[
\bigl| P_M^{\bo,x}(C_i \times C_j)-P_M^{\mathbf o}(C_i)P_M^{\mathbf o}(C_j)
\bigr|+ \bigl|P_M^{\bo,x}(C_j \times
\mathbb{M})-P_M^{\mathbf o}(C_j) \bigr| +
\bigl|P_M^{\bo,x}(C_i \times\mathbb{M})-P_M^{\mathbf o}(C_i)\bigr|,
\]
which, by (\ref{intconalf}), is integrable w.r.t. $\alpha_{\X,\mathrm{red}}^{(2)} $.
Hence, (\ref{dartautil}) follows by (\ref{ineqHeinrichPawlas}) and
Lebesgue's
dominated convergence theorem.

\subsection{Proof of Theorem \texorpdfstring{\protect\ref{theconvar}}{3.4}}
We again expand the parentheses in the second term of the estimator
$(\widehat{\sigma}_{ij}^{(1)})_k$ defined by (\ref{esttauone}) and
express the expectations in terms of $P_M^{\bo,y}$ and
$\alpha_{\X,\mathrm{red}}^{(2)}$. Using the obvious relation $\gamma
_k(y)=\int_{\R^d} \ind_{W_k}(x) \ind_{W_k}(y+x) \,\rmd x$ we find that,
for any $C_i, C_j \in\calB(\mathbb{M}) $,
\begin{eqnarray*}
&&
\E\sum^{\neq} _{p,q \geq1}\frac{
\ind_{W_k}(X_p)\ind_{W_k}(X_q)\ind_{C_i}(M_p)\ind
_{C_j}(M_q)}{|(W_k-X_p)\cap
(W_k-X_q)|} \\
&&\quad=
\int_{(\R^d)^2} \frac{\ind_{W_k}(x) \ind_{W_k}(y)
P_M^{x,y}(C_i \times C_j)}{\gamma_k(y-x)} \alpha_\X^{(2)}
\bigl(\rmd (x,y)\bigr)
\\
&&\quad= \lambda\int_{\R^d} \frac{P_M^{\bo,y}(C_i \times
C_j)}{\gamma_k(y)} \int
_{\R^d} \ind_{W_k}(x) \ind_{W_k}(y+x) \,\rmd x \alpha_{\X,\mathrm{red}}^{(2)}(\rmd y) \\
&&\quad=\lambda\int
_{\R^d}P_M^{\bo,y}(C_i
\times C_j) \alpha_{\X,\mathrm{red}}^{(2)}(\rmd y).
\end{eqnarray*}
As in the proof of Theorem \ref{thereptau} after summarizing all
terms we obtain that
\begin{eqnarray*}
\E\bigl(\widehat{\sigma}_{ij}^{(1)}\bigr)_k&=&
\lambda\bigl(P_M^{\mathbf o}(C_i \cap
C_j)-P_M^{\mathbf o}(C_i)P_M^{\mathbf
o}(C_j)
\bigr)
\\
&&{}+\lambda\int_{\R^d} \bigl(P_M^{\bo,x}(C_i
\times C_j)- P_M^{\bo,x}(C_i
\times\mathbb{M})P_M^{\mathbf o}(C_j)
\\
&&\hspace*{36pt}{}- P_M^{\bo,x}(C_j\times\mathbb{M})
P_M^{\mathbf o}(C_i) + P_M^{\mathbf o}(C_i)
P_M^{\mathbf o}(C_j) \bigr) \alpha^{(2)}_{\X,\mathrm{red}}(\rmd x),
\end{eqnarray*}
which by comparison to (\ref{dartautil}) yields that $\E(\widehat
{\sigma}_{ij}^{(1)})_k=
\sigma_{ij} $.
The asymptotic unbiasedness of $(\widehat{\sigma}_{ij}^{(2)})_k$ is
rapidly seen by (\ref
{covmattau}) and the equality
$\E(\widehat{\sigma}_{ij}^{(2)})_k= \operatorname{Cov} (
Y_k(C_i),Y_k(C_j) )= \E
Y_k(C_i)Y_k(C_j) $, which follows
directly from (\ref{streu2}).

\subsection{Proof of Theorem \texorpdfstring{\protect\ref{thmestTauThree}}{3.5}} \label{secproofCovEst}

Since $\E( \sigma_{ij} - (\widehat{\sigma}_{ij}^{(3)})_k
)^2 = \Var(\widehat{\sigma}_{ij}^{(3)})_k+
( \sigma_{ij} - \E(\widehat{\sigma}_{ij}^{(3)})_k )^2$
we have to show that
%
\begin{equation}
\label{formVarE} \E\bigl(\widehat{\sigma}_{ij}^{(3)}
\bigr)_k\mathop{\longrightarrow} _{k \to\infty}
\sigma_{ij} \quad\mbox{and}\quad \Var\bigl(\widehat{\sigma}_{ij}^{(3)}
\bigr)_k\mathop{\longrightarrow} _{k \to\infty}0.
\end{equation}
For notational ease, we put
\begin{eqnarray*}
m(u,v) &=& \bigl(\ind_{C_i}(u) - P_M^{\mathbf o}(C_i)
\bigr) \bigl(\ind_{C_j}(v) - P_M^{\mathbf
o}(C_j)
\bigr),\qquad a_k = b_k |W_k|^{1/d},
\\
r_k(x,y) &=& \frac{\ind_{W_k}(x) \ind_{W_k}(y)}{\gamma_k(y-x)} w \biggl
(\frac{\|y-x\|}{a_k} \biggr)
\quad\mbox{and}\quad \tau_k = \sum_{p,q\ge1}^{\neq}
r_k(X_p, X_q) m(M_p,
M_q).
\end{eqnarray*}
Hence, together with (\ref{ergempmark}) and (\ref{ypskahtil}) we
may rewrite $(\widehat{\sigma}_{ij}^{(3)})_k$ as follows:
%
\begin{equation}
\label{streu3} \bigl(\widehat{\sigma}_{ij}^{(3)}
\bigr)_k=\frac{1}{\sqrt{|W_k|}} Y_k(C_i\cap
C_j)+ \widehat{\lambda}_k \bigl(P_M^{\mathbf o}(C_i
\cap C_j)-P_M^{\mathbf o}(C_i)P_M^{\mathbf
o}(C_j)
\bigr)+\tau_k.
\end{equation}
Using the definitions and relations
(\ref{defalfisk})--(\ref{condRadonNikodym}) and
$\int_{\R^d}r_k(x,y+x)\,\rmd x = w (\|y\|/a_k )$ we find that the
expectation $\E\tau_k$ can be expressed by
\begin{eqnarray*}
&&\int_{(\R^d\times{\mathbb M})^2}r_k(x,y) m(u,v)
\alpha_{\X_M}^{(2)} \bigl(\rmd (x,u,y,v) \bigr)
\\
&&\quad= \lambda\int_{\R^d} \int_{{{\mathbb M}}^2}
m(u,v)P_M^{\bo,y} \bigl(\rmd (u,v) \bigr) w \biggl(
\frac{\|y\|}{a_k} \biggr) \alpha_{\X,\mathrm{red}}^{(2)} (\rmd y ).
\end{eqnarray*}
The inner integral $\int_{{\mathbb M}^2} m(u,v) P_M^{\bo,y}
(\rmd (u,v) )$ coincides with the integrand occurring in (\ref
{dartautil}) and this term is integrable w.r.t. $\alpha_{\X,\mathrm{red}}^{(2)}$
due to
(\ref{intconalf}) which in turn is a consequence of (\ref
{condConsistencyTau}) and
Lemma \ref{lemintgam}. Hence, by Condition \hyperref[condWB]{$(wb)$} and the dominated
convergence theorem,
we arrive at
\begin{eqnarray*}
\E\tau_k &\displaystyle \mathop{\longrightarrow} _{k \to\infty}& \lambda\int
_{\R^d}\int_{{{\mathbb M}}^2}m(u,v)
P_M^{\bo,y} \bigl(\rmd (u,v) \bigr) \alpha_{\X,\mathrm{red}}^{(2)}
(\rmd y )
\\
&=& \sigma_{ij} - \lambda\bigl( P_M^{\mathbf o}(C_i
\cap C_j)- P_M^{\mathbf o}(C_i)P_M^{\mathbf o}
( C_j) \bigr).
\end{eqnarray*}
The definitions of $\widehat{\lambda}_k$ and $Y_k(\cdot)$ by (\ref
{ergempmark}) and
(\ref{ypskahtil}), respectively, reveal that
$\E\widehat{\lambda}_k = \lambda$ and $\E Y_k(C_i\cap C_j) = 0
$. This
combined with the last limit and (\ref{streu3}) proves the first relation
of (\ref{formVarE}).
To verify the second part of (\ref{formVarE}) we apply the Minkowski
inequality to the rhs of (\ref{streu3}) which yields the estimate
\[
\bigl(\Var\bigl(\widehat{\sigma}_{ij}^{(3)}
\bigr)_k \bigr)^{1/2} \le|W_k|^{-1/2}
\bigl(\Var Y_k(C_i\cap C_j)
\bigr)^{1/2} + (\Var\widehat{\lambda}_k )^{1/2} +
(\Var\tau_k )^{1/2}.
\]
The first summand on the rhs tends to 0 as $k \to\infty$ since $\E Y_k(C)^2$
has a finite limit for any $C\in{\mathcal B}({\mathbb M})$ as shown in
Theorem \ref{thereptau}
under condition (\ref{intconalf}). The second summand is easily seen to
disappear as $k \to\infty$ if (\ref{intcummea}) is fulfilled, see,
for example, \cite{Hei94,Heinrich08} or \cite{Hei10}.
Condition (\ref{condConsistencyTau}) implies both (\ref
{intcummea}) and
(\ref{intconalf}), see Lemma \ref{lemintgam}. Therefore, it
remains to show that
$\Var\tau_k \longrightarrow0$ as $k \to\infty$. For this
purpose, we employ the variance formula (\ref{streuung}) stated in
Lemma \ref{LemmaFormelMonster} in the special case $f(x,y,u,v) =
r_k(x,y) m(u,v) $. In this way, we get the decomposition
$\Var\tau_k = I_k^{(1)} + I_k^{(2)} + I_k^{(3)} $, where
$I_k^{(1)}$, $I_k^{(2)}$ and
$I_k^{(3)}$ denote the three multiple integrals on the rhs of (\ref
{streuung}) with $f(x,y,u,v)$
replaced by the product $r_k(x,y) m(u,v) $. We will see that the integrals
$I_k^{(1)}$ and $I_k^{(2)}$ are easy to estimate only by using (\ref
{intcummea}) and (\ref{intconalf}) while in order to show that
$I_k^{(3)}$ tends to $0$ as $k \to\infty$, the full strength of the
mixing condition (\ref{condConsistencyTau}) must be exhausted. Among
others we use repeatedly the estimate
%
\begin{equation}\label{formgammak}
\frac{1}{\gamma_k(a_k y)}\le\frac{2}{|W_k|} \qquad\mbox{for } y \in B(\bo,
r_w),
\end{equation}
which follows directly from (\ref{ineqHeinrichPawlas}) and the choice of
$\{b_k\}$ in (\ref{condbk1}). The definition of $I_k^{(1)}$ together with
(\ref{formgammak}) and $\alpha_{\X, \mathrm{red}}^{(2)}(\rmd x)= \gamma
_{\X, \mathrm{red}}^{(2)}(\rmd x)+\lambda\,\rmd x$ yields
\begin{eqnarray*}
\bigl|I_k^{(1)}\bigr| &\le& 2 \int_{(\R^d)^2} \bigl(r_k(x_1,x_2) \bigr)^2
\alpha_\X^{(2)} \bigl(\rmd (x_1,x_2) \bigr) = 2 \lambda\int_{\R^d}
\frac{1}{\gamma_k(y)} w^2 \biggl(\frac{\|y\|}{a_k}
\biggr)\alpha_{\X,\mathrm{red}}^{(2)}(\rmd y)
\\
&\le& \frac{4 \lambda}{|W_k|} \biggl( m_w^2 \bigl|
\gamma_{\X,\mathrm{red}}^{(2)}\bigr|\bigl(\R^d\bigr) + \lambda
a_k^d \int_{\R^d} w^2\bigl(
\|y\|\bigr)\,\rmd y \biggr) \mathop{\longrightarrow} _{k \to\infty}0,
\end{eqnarray*}
where the convergence results from Condition \hyperref[condWB]{$(wb)$} and (\ref
{condConsistencyTau}), which implies
$|\gamma_{\X, \mathrm{red}}^{(2)}|(\R^d)<\nolinebreak\infty$ by virtue of
Lemma \ref{lemintgam}. Analogously, using besides (\ref
{formgammak}) and Condition \hyperref[condWB]{$(wb)$} the relations
\[
w \biggl(\frac{\|x\|}{a_k} \biggr) \le m_w \ind_{[-\lceil a_k
r_w\rceil,
\lceil a_k r_w\rceil]^d}(x)
\quad\mbox{and}\quad W_k \subseteq\bigcup_{z\in{\overline H_k}}E_z
\]
with the notation introduced in Section \ref{secPalmMarkDistr}
we obtain that
\begin{eqnarray*}
\bigl|I_k^{(2)}\bigr| &\le& 4 \int_{(\R^d)^3}
r_k(x_1,x_2) r_k(x_1,
x_3) \alpha_\X^{(3)} \bigl(\rmd (x_1,x_2,x_3) \bigr)
\\
&\le& \frac{16 m_w^2}{|W_k|^2}\sum_{z \in{\overline H_k}} \alpha
_\X^{(3)} \bigl(\bigl(E_z\oplus\bigl[-\lceil
a_k r_w\rceil, \lceil a_k
r_w\rceil\bigr]^d\bigr) \times\bigl(E_z
\oplus\bigl[-\lceil a_k r_w\rceil, \lceil
a_k r_w\rceil\bigr]^d\bigr) \times
E_z \bigr).
\end{eqnarray*}
Since the cube $E_z\oplus[-\lceil a_k r_w\rceil, \lceil a_k r_w\rceil]^d$
decomposes into $(2\lceil a_k r_w\rceil+1)^d$ disjoint unit cubes and
$\alpha_\X^{(3)}(E_{z_1} \times E_{z_2} \times E_{z_3} )\le\E(\X
(E_{\bo}))^3$ by H\"older's inequality, we may proceed with
\[
\bigl|I_k^{(2)}\bigr| \le\frac{16 m_w^2}{|W_k|^2} \#{\overline
H}_k \bigl(2\lceil a_k r_w\rceil+1
\bigr)^{2d} \E\bigl(\X(E_{\bo})\bigr)^3 \le
c_1 b_k^{2d} |W_k| \mathop{
\longrightarrow} _{k \to\infty}0.
\]
Here we have used the moment condition in (\ref{condConsistencyTau}),
(\ref{CASgridConvergence}), 
and the assumptions (\ref{condbk1}) imposed on the sequence $\{b_k\}
$.

In order to prove that $I_k^{(3)}$ vanishes as $k \to\infty$, we first
evaluate the inner integrals over the product $m(u_1,u_2) m(u_3,u_4)$
with $m(u,v)= (\ind_{C_i}(u)-P^{\bo}_M(C_i) ) (\ind_{C_j}(v)
-P^{\bo}_M(C_j) )$ so that $I_k^{(3)}$ can be written as linear
combination
of 16 integrals taking the form
\begin{eqnarray*}
J_k &=& \int_{(\R^d)^2}\int_{(\R^d)^2}
r_k(x_1,x_2) r_k(x_3,x_4)\\
&&\hspace*{48.4pt}{}\times
\Biggl[ P^{x_1,x_2, x_3, x_4}_M\Biggl(\bigtimes_{r=1}^4D_r
\Biggr) \alpha_\X^{(4)} \bigl(\rmd (x_1,x_2,
x_3, x_4) \bigr)
\\
&&\hspace*{65.5pt}{} - P^{x_1,x_2}_M(D_1\times D_2)
P^{x_3,x_4}_M(D_3\times D_4)
\alpha_\X^{(2)} \bigl(\rmd (x_1,x_2)
\bigr) \alpha_\X^{(2)} \bigl(\rmd (x_3,x_4)
\bigr) \Biggr]
\\
&=& \int_{\timesop_{r=1}^4(\R^d\times D_r)} r_k(x_1,x_2)
r_k(x_3,x_4) \bigl( \alpha_{\X_M}^{(4)}
- \alpha_{\X_M}^{(2)}\times\alpha_{\X
_M}^{(2)}
\bigr) \bigl(\rmd (x_1,u_1,\ldots,x_4,u_4)
\bigr),
\end{eqnarray*}
where the mark sets $D_1, D_3 \in\{C_i, {\mathbb M}\}$ and $D_2, D_4
\in\{C_j, {\mathbb M}\}$ are fixed in what follows and the signed measure
$\alpha_{\X_M}^{(4)} - \alpha_{\X_M}^{(2)}\times\alpha_{\X
_M}^{(2)}$ on
${\mathcal B}((\R^d\times{\mathbb M})^4)$ (and its total variation
measure $ |\alpha_{\X_M}^{(4)} - \alpha_{\X_M}^{(2)}\times
\alpha_{\X_M}^{(2)} | $) come into play by virtue of the definition
(\ref{condRadonNikodym}) for the $m$-point Palm mark distribution in
case $m=2$
and $m=4 $.

As $|z_1-z_2| > \lceil a_k r_w \rceil$ (where, as above, $|z|$ denotes
the maximum norm of $z\in{\mathbb Z}^d$) implies $\|x_2-x_1\| > a_k
r_w$ and thus $r_k(x_1,x_2) = 0$ for all
$x_1 \in E_{z_1}, x_2 \in E_{z_2}$, we deduce from (\ref
{formgammak}) together with Condition \hyperref[condWB]{$(wb)$} and the abbreviation
$N(a_k) = (1+c_0)(\lceil a_k r_w \rceil+ 1)$
(where $c_0$ is from (\ref{betstaone})) that
%
\begin{equation}\label{formIk3Sum1}
|J_k|\le\frac{4 m_w^2}{|W_k|^2} \Biggl(\sum
_{n=0}^{\lceil N(a_k)
\rceil}+ \sum_{n > \lceil N(a_k) \rceil}
\Biggr)\sum_{(z_1,z_2)\in S_k} \sum_{(z_3,z_4)\in S_{k,n}(z_1)}
V_{z_1, z_2, z_3, z_4},
\end{equation}
where $S_k=\{(u,v) \in{\overline H_k}\times{\overline H_k}\dvt  |u-v|\le
\lceil a_k r_w \rceil\}, S_{k,n}(z)=\{(z_1,z_2)\in S_k\dvt
\min_{i=1,2}|z_i-z|=n\}$ and $V_{z_1, z_2, z_3, z_4} =
|\alpha_{\X_M}^{(4)} - \alpha_{\X_M}^{(2)}\times\alpha_{\X
_M}^{(2)} |
(\bigtimes_{r=1}^4 (E_{z_r}\times D_r) )$ for any
$z_1,\ldots,z_4 \in{\mathbb Z}^d $.
Obviously, for any fixed $z \in{\overline H_k}$, at most
$2 (\lceil N(a_k) \rceil+1)^d (2 \lceil N(a_k) \rceil+1)^d$ pairs
$(z_3,z_4)$ belong to $\bigcup_{n=0}^{\lceil N(a_k) \rceil} S_{k,n}(z)$
and the number of pairs $(z_1,z_2)$ in $S_k$ does not exceed the
product $\#{\overline H_k} (2 \lceil a_k r_w \rceil+1)^d$. Finally,
remembering that $a_k = b_k |W_k|^{1/d}$ and using
the evident estimate $V_{z_1, z_2, z_3, z_4} \le2 \E(\X(E_{\bo}))^4$
together with (\ref{CASgridConvergence}) and Condition \hyperref[condWB]{$(wb)$}, we
arrive at
\[
\frac{4 m_w^2}{|W_k|^2} \sum_{(z_1,z_2)\in S_k} \sum
_{n=0}^{\lceil N(a_k) \rceil}\sum_{(z_3,z_4)\in S_{k,n}(z_1)}
V_{z_1, z_2, z_3, z_4} \le c_2 \frac{\#{\overline H_k}}{|W_k|^2} \bigl(b_k^d
|W_k| \bigr)^3 \mathop{\longrightarrow}
_{k \to
\infty}0.
\]
It remains to estimate the sums on the rhs of (\ref{formIk3Sum1})
running over $n > \lceil N(a_k) \rceil$.
For the signed measure $\alpha^{(4)}_{\X_M} - \alpha^{(2)}_{\X
_M}\times\alpha^{(2)}_{\X_M}$ we consider the Hahn decomposition
$H^+, H^- \in{\mathcal B}((\R^d \times{\mathbb M})^4)$ yielding
positive (negative) values on subsets of $H^+$($H^-$). Recall that
$K_a=[-a,a]^d$.
For fixed $z_1 \in{\overline H_k}$, $z_2\in{\overline H_k} \cap
(K_{\lceil a_k r_w \rceil}+z_1)$
and $(z_3,z_4)\in S_{k,n}(z_1)$, we now consider the decompsition
$V_{z_1, z_2, z_3, z_4} =
V^{+}_{z_1, z_2, z_3, z_4} + V^{-}_{z_1, z_2, z_3, z_4}$ with
\[
V^{\pm}_{z_1, z_2, z_3, z_4} = \pm\bigl(\alpha_{\X_M}^{(4)}
- \alpha_{\X_M}^{(2)}\times\alpha_{\X_M}^{(2)}
\bigr) \Biggl(H^{\pm} \cap\bigtimes_{r=1}^4
(E_{z_r}\times D_r) \Biggr).
\]
Since $(z_3,z_4)\in S_{k,n}(z_1)$ means that $z_3 \in{\overline H_k}
\cap(K^c_n+z_1 )$, where $K^c_a=\R^d \setminus K_a
$, and $z_4\in{\overline H_k}\cap(K_{\lceil a_k r_w \rceil
}+z_3 )\cap(K^c_n+z_1 )$, we define MPPs $Y_k$
and $Y_n'$ as the restrictions of $\X_M$ to
$(K_{\lceil a_k r_w \rceil+1/2}+z_1)\times{\mathbb M}$ and
$(K^c_{n-1/2}+z_1)\times{\mathbb M} $, respectively. Let
furthermore $\widetilde Y_k$ and $\widetilde Y_n'$ be copies of $Y_k$
and $Y_n'$ which are independent.
Next, we define functions $f^{+}(Y_k,Y_n^{\prime})$ and
$f^{-}(Y_k,Y_n^{\prime})$ by
\[
f^{\pm}\bigl(Y_k,Y_n^{\prime}\bigr) =
\sum^{\neq}_{p,q\ge1}\sum
^{\neq
}_{s,t\ge1} \ind_{\pm}\bigl(X_p,M_p,X_q,M_q,X_s',M_s',X_t',M_t'
\bigr),
\]
where $\ind_{\pm}(\cdots)$ denote the indicator functions of the sets
$H^{\pm} \cap\bigtimes_{r=1}^4 (E_{z_r}\times D_r)$ so that we get
\[
V^{\pm}_{z_1, z_2, z_3, z_4} = \E f^{\pm}\bigl(Y_k,Y_n^{\prime}
\bigr)- \E f^{\pm}\bigl(\widetilde Y_k, \widetilde
Y_n'\bigr) \qquad\mbox{for } (z_1,z_2)
\in S_k, (z_3,z_4)\in S_{k,n}(z_1).
\]
Hence, having in mind the stationarity of $\X_M$, we are in a position to
apply the covariance inequality (\ref{covinebet}), 
which provides for $\eta> 0$ and $n > \lceil N(a_k) \rceil$ that
%
\begin{eqnarray}\label{Vestimate}
V^{\pm}_{z_1, z_2, z_3, z_4} &\le& 2 \bigl( \beta\bigl(
\calA(K_{\lceil a_k r_w \rceil+1/2}+z_1), \calA\bigl(K^c_{n-1/2}+z_1
\bigr) \bigr) \bigr)^{{\eta}/({1+\eta
})}
\nonumber
\\
&&{}\times \Biggl( \E\Biggl(\prod_{r=1}^2
\X_M(E_{z_r}\times D_r)
\Biggr)^{2+2\eta} \E\Biggl(\prod_{r=3}^4
\X_M(E_{z_r}\times D_r)
\Biggr)^{2+2\eta} \Biggr)^{{1}/({2+2\eta})}\quad
\\
&\le& 2 \bigl(\beta^*_{\X_M}\bigl(n-\lceil a_kr_w
\rceil-1\bigr) \bigr)^{{\eta}/({1+\eta})} \bigl(\E\X(E_{\bo})^{4+4\eta
}
\bigr)^{1/(1+\eta)}.\nonumber
\end{eqnarray}
In the last step, we have used the Cauchy--Schwarz inequality and the definition
of the $\beta$-mixing rate $\beta^*_{\X_M}$ together with constant
$c_0$ in (\ref{betstaone}). Finally,
setting $\eta= \delta/4$ with $\delta> 0$ from (\ref{condConsistencyTau})
the estimate (\ref{Vestimate}) enables us to derive the following
bound of
that part on the rhs of (\ref{formIk3Sum1}) connected with the series over
$n > \lceil N(a_k) \rceil$:
\[
c_3 \frac{\#{\overline H_k}}{|W_k|^2} \bigl(2\lceil a_kr_w
\rceil+1\bigr)^{2d} \sum_{n > \lceil N(a_k) \rceil} \bigl(
(2n+1)^d - (2n-1)^d \bigr) \bigl(\beta^*_{\X_M}
\bigl(n-\lceil a_kr_w\rceil-1\bigr)
\bigr)^{{\delta
}/({4+\delta})}.
\]
Combining $a_k = b_k|W_k|^{1/d}$ and (\ref{CASgridConvergence})
with condition (\ref{condConsistencyTau}) and the choice of $\{b_k\}$
in (\ref{condbk1}), it is
easily checked that the latter expression and thus $J_k$ tend to 0 as
$k \to\infty$. This completes the proof of Theorem \ref{thmestTauThree}.

\section{Examples} \label{secexamples}
\subsection{$m$-dependent marked point processes}
A stationary MPP $\X_M$ is called $m$-\emph{dependent} if, for any
$B, B'\in\mathcal{B}(\R^d)$, the $\sigma$-algebras $\mathcal
{A}_{\X_M}(B)$
and $\mathcal{A}_{\X_M}(B')$ are stochastically
independent if $\inf\{|x-y|\dvt  x \in B, y\in B'\} > m$ or, equivalently,
\[
\beta\bigl(\mathcal{A}_{\X_M}(K_a),
\mathcal{A}_{\X
_M}\bigl(K^c_{a+b}\bigr) \bigr) = 0
\qquad\mbox{for } b > m \mbox{ and } a > 0.
\]
In terms\vspace*{1pt} of the corresponding mixing rates this means that $\beta_{\X
_M}^*(r)=\beta_{\X_M}^{**}(r)=0$ if $r > m $. For $m$-dependent MPPs
$\X_M$, it is evident that Condition
\hyperref[condBeDe]{${\beta(\delta)}$} in Theorem \ref{theasynor} is
only meaningful for $\delta= 0 $, that is, $\E\X([0,1]^d)^2 < \infty$.
This condition also implies (\ref{intcummea}) and (\ref{intconalf}).
Likewise, the assumption (\ref{condConsistencyTau}) of Theorem
\ref{thmestTauThree} reduces to $\E\X([0,1]^d)^4 < \infty$ which
suffices to prove the $L^2$-consistency of the empirical covariance
matrix $\widehat{\bolds\Sigma}_k^{(3)} $.

\subsection{Geostatistically marked point processes}
Let $\X= \sum_{n\ge1}\delta_{X_n}$ be an unmarked simple PP on $\R
^d$ and
$M = \{M(x), x \in\R^d\}$ be a measurable random field on $\R^d$ taking
values in the Polish mark space ${\mathbb M}$. Further assume that $\X
$ and $M$
are stochastically independent over a common probability space
$(\Omega, {\mathcal A}, {\mathbb P})$. An MPP $\X_M = \sum_{n\ge
1}\delta_{(X_n,M_n)}$ with atoms $X_n$ of $\X$ and marks $M_n =
M(X_n)$ is called \emph{geostatistically marked}. Equivalently, the
random counting measure
$\X_M \in{\mathsf N}_{\mathbb M}$ can be represented by means of the Borel
sets $M^{-1}(C) = \{x\in\R^d\dvt  M(x) \in C\}$ (if $C\in{\mathcal
B}({\mathbb M})$) by
%
\begin{equation}
\label{geostat} \X_M(B \times C) = \X\bigl(B \cap
M^{-1}(C)\bigr) \qquad\mbox{for } B \times C \in{\mathcal B}\bigl(
\R^d\bigr)\times{\mathcal B}({\mathbb M}).
\end{equation}
Obviously, if both the PP $\X$ and the mark field $M$ are stationary
then so is
$\X_M$ and vice versa. Furthermore, the $m$-dimensional distributions of
$M$ coincide with the $m$-point Palm mark distributions of $\X_M$. The
following lemma allows to estimate the $\beta$-mixing
coefficient (\ref{defbetemm}) by the sum of the
corresponding coefficients of the PP $\X$ and the mark field $M$.

\begin{Lemma} Let the MPP $\X_M$ be defined by (\ref{geostat}) with
an unmarked
PP and a random mark field $M$ being stochastically independent of each other.
Then, for any $B, B'\in\mathcal{B}(\R^d) $,
%
\begin{equation}
\label{betmixcoef} \beta\bigl(\calA_{\X_M}(B),\calA_{\X_M}
\bigl(B^\prime\bigr) \bigr) \le\beta\bigl(\calA_\X(B),
\calA_\X\bigl(B^\prime\bigr) \bigr) + \beta\bigl(
\calA_M(B),\calA_M\bigl(B^\prime\bigr) \bigr),
\end{equation}
where the $\sigma$-algebras $\calA_\X(B),\calA_\X(B')$ and $\calA
_M(B),\calA_M(B')$ are
generated by the restriction of $\X$ and $M$, respectively, to the
sets $B,B'$.
\end{Lemma}
To sketch a proof for (\ref{betmixcoef}), we regard the differences
$\Delta(A_i,A_j')={\mathbb P}(A_i\cap A_j')- {\mathbb P}(A_i)
{\mathbb P}(A_j')$ for two finite partitions $\{A_i\}$ and $\{A'_j\}$
of $\Omega$ consisting of events of the form
\begin{eqnarray*}
A_i &=& \bigcap_{p=1}^{k}
\bigl\{\X_M(B_p\times C_p)\in
\Gamma_{p,i}\bigr\},\\
A_j' &=& \bigcap
_{q=1}^{{\ell}}\bigl\{\X_M
\bigl(B_q'\times C_q'\bigr)\in
\Gamma'_{q,j}\bigr\} \qquad\mbox{with } \Gamma_{p,i},
\Gamma'_{q,j} \subseteq{\mathbb Z}_+
\end{eqnarray*}
with pairwise disjoint bounded Borel sets $B_1,\ldots,B_k \subseteq B$ and
$B_1',\ldots,B_{\ell}' \subseteq B'$. This suffices since the
supremum in (\ref{defbetemm}) does not change if the sets $A_i$ and
$A_j^{\prime}$ belong to semi-algebras generating $\calA_{\X_M}(B)$ and
$\calA_{\X_M}(B^\prime)$, respectively. Making use of
(\ref{geostat}) combined with the independence assumption yields the identity
\begin{eqnarray*}
\Delta\bigl(A_i,A_j'\bigr) &=& \int
_{\Omega}\int_{\Omega} ({\mathbb
P}_{\calA
_\X(B)\otimes\calA_\X(B^\prime)}-{\mathbb P}_{\calA_\X(B)}\times
{\mathbb P}_{\calA_\X(B^\prime)}
) \bigl(A_i\cap A_j'\bigr) \,\rmd  {
\mathbb P}_{\calA_M(B)\otimes\calA_M(B^\prime)}
\\
&&{}+ \int_{\Omega}\int_{\Omega}{\mathbb
P}_{\calA
_\X(B)}(A_i) {\mathbb P}_{\calA_\X(B^\prime)}
\bigl(A_j'\bigr) \,\rmd  ({\mathbb P}_{\calA_M(B)\otimes\calA_M(B^\prime)}-
{\mathbb P}_{\calA
_M(B)}\times{\mathbb P}_{\calA_M(B^\prime)} ),
\end{eqnarray*}
which by (\ref{defbetemm}) and the integral form of the total
variation confirms (\ref{betmixcoef}).

\subsection{Cox processes on the boundary of germ-grain models}\label
{seccoxBoolean}
Let $\Xi= \bigcup_{n \ge1} (\Xi_n+Y_n)$ be a \emph{germ-grain
model}, see, for example, \cite{Hei99}, governed by some stationary
unmarked PP
${\mathbf Y}=\sum_{n\ge1} \delta_{Y_n}$ in $\R^d$ with intensity
$\lambda>0$ and a sequence $\{\Xi_n\}_{n\ge1}$ of independent copies
of some
random convex, compact set $\Xi_0$ (such that ${\mathbb P}(\bo\in\Xi_0)=1$)
called \emph{typical grain}. With the radius functional $\|\Xi_0\|=
\sup\{
\|x\|\dvt  x \in\Xi_0\}$, the condition $\E\|\Xi_0\|^d < \infty$
ensures that $\Xi$ is a random closed set. The germ-grain model is
called \emph{Boolean model} if the PP ${\mathbf Y}$ is Poisson.
We consider a marked Cox process $\X_M$,
where the unmarked Cox process $\X=\sum_{n\ge1}\delta_{X_n}$ is
concentrated on the boundary $\partial\Xi$ of $\Xi$ with
random intensity measure being proportional to the $(d-1)$-dimensional
Hausdorff measure $\mathcal{H}_{d-1}$ on $\partial\Xi$.
As marks $M_n$ we take the outer unit normal vectors at the points
$X_n \in\partial\Xi$, which are (a.s.) well defined for $n\ge1$ due
to the assumed convexity of $\Xi_0$.
This example with marks given by the orientation of outer normals in
random boundary points may occur rather specific. However, in this way
our asymptotic results may be used to construct asymptotic tests for
the fit of a Boolean model to a given dataset w.r.t. its rose of
directions. For instance, if the typical grain
$\Xi_0$ is rotation-invariant (implying the isotropy of $\Xi$), then
the Palm mark distribution $P_M^{\mathbf o}$ of the stationary MPP $\X
_M = \sum_{n\ge1}\delta_{(X_n,M_n)}$ is the uniform distribution on the unit
sphere ${\mathbb S}^{d-1}$ in $\R^d$. We will now discuss assumptions
ensuring that Condition \hyperref[condBeDe]{${\beta(\delta)}$} and (\ref
{condConsistencyTau}) hold, which are required for our CLT (\ref
{asyequnor}) and the consistent estimation of the covariances (\ref
{covmattau}), respectively.
Using Lemmas 5.1 and 5.2 in \cite{Hei99} (with improved constants), we
obtain that
\begin{eqnarray*}
&&
\beta\bigl(\calA_{\X_M}(K_a),\calA_{\X_M}
\bigl(K^c_{a+b}\bigr) \bigr) \\
&&\quad\le\beta\bigl(
\calA_{\mathbf Y}(K_{a+b/4}),\calA_{\mathbf
Y}
\bigl(K^c_{a+3 b/4}\bigr) \bigr)
\\
&&\qquad{} + \lambda2^{d+1} \biggl( \biggl(1+\frac{4a}{b}
\biggr)^{d-1} + \biggl(3+\frac{4a}{b} \biggr)^{d-1}
\biggr) \E\|\Xi_0\| ^d\ind\biggl\{\|\Xi_0
\|\ge\frac{b}{4}\biggr\}
\end{eqnarray*}
for $a,b \ge1/2$. According to (\ref{betstaone}) with $c_0=4$, we
may thus define the $\beta$-mixing rates
$\beta_{\X_M}^*(r)$ and $\beta_{\X_M}^{**}(r)$ for $r \ge2$ to be
\begin{eqnarray*}
\beta_{\X_M}^*(r) &=& \beta_{\mathbf Y}^*\biggl(\frac{r}{2}
\biggr) + c_4 \E\|\Xi_0\|^d\ind\biggl\{\|
\Xi_0\|\ge\frac{r}{4}\biggr\} \ge\sup_{a \in[1/2,
r/4]}
\beta\bigl(\calA_{\X_M}(K_a),\calA_{\X_M}
\bigl(K^c_{a+r}\bigr) \bigr),
\\
\beta_{\X_M}^{**}(r) &=& 2^{d-1} \beta_{\mathbf Y}^{**}
\biggl(\frac
{r}{2}\biggr) + c_4 \frac{4^{d-1}}{r^{d-1}} \E\|
\Xi_0\|^d\ind\biggl\{\|\Xi_0\|\ge
\frac{r}{4}\biggr\} \ge\sup_{a \ge
r/4}\frac
{\beta(\calA_{\X_M}(K_a),\calA_{\X_M}(K^c_{a+r}) )}{a^{d-1}}
\end{eqnarray*}
with $c_4 =\lambda4^d (1 + 2^{d-1})$ and rate functions $\beta
_{\mathbf Y}^{*}(r), \beta_{\mathbf Y}^{**}(r)$ which are defined in
analogy to (\ref{betstaone}) for $c_0=4$.

It is easily seen that
\[
\E\|\Xi_0\|^{2d} < \infty\quad\mbox{and}\quad {(A):}\quad
r^{2d-1}\beta_{\mathbf
Y}^{**}(r) \lonr0
\]
imply $r^{2d-1}\beta_{\X_M}^{**}(r) \lonr0$. Moreover,
\[
{(B_{\delta,p}):}\quad \E\|\Xi_0\|^{2d(p+\delta)/\delta} < \infty
\quad\mbox{and}\quad {(C_{\delta,p}):}\quad \int_1^\infty
r^{d-1} \bigl( \beta^*_{{\mathbf Y}}(r) \bigr)^{\delta/(2p+\delta)}
\,\rmd r < \infty
\]
ensure
$\int_1^\infty r^{d-1} ( \beta^*_{\X_M}(r)
)^{\delta/(2p+\delta)} \,\rmd r < \infty$ for any $p \ge0$
and $\delta> 0$.
Further, the random intensity measure of $\X$ on $E_{\bo}$ and thus
also $\X(E_{\bo})$
has moments of order $q \ge1$ if $\E{\mathbf Y}(E_{\bo})^q < \infty
$ and $\E\|\Xi_0\|^d < \infty$.
Now we are in a position to express Condition \hyperref[condBeDe]{${\beta(\delta)}$} and
(\ref{condConsistencyTau}) by conditions on $\Xi_0$ and ${\mathbf Y}$.

\begin{Lemma}\label{lemcondcox} For the above-defined stationary
marked Cox process $\X_M$ on the
boundary of the germ-grain model $\Xi$ generated by the PP ${\mathbf
Y}$ and typical grain $\Xi_0$, the assumptions of
Theorem~\ref{theasynor}, respectively, Theorem \ref{thmestTauThree} are
satisfied whenever, for some $\delta> 0$,
\[
\E{\mathbf Y}(E_{\bo})^{2+\delta} < \infty, {(A)},
{(B_{\delta,1})}, {(C_{\delta,1})} \mbox{, respectively, } \E{\mathbf
Y}(E_{\bo})^{4+\delta} < \infty, {(B_{\delta,2})},
{(C_{\delta,2})}.
\]
\end{Lemma}
\begin{Remark*}
If the stationary PP ${\mathbf Y}$ of germs is Poisson
the conditions $\E{\mathbf Y}(E_{\bo})^{4+\delta} < \infty $, ${(A)}$
and ${(C_{\delta,2})}$ are trivially satisfied for any $\delta> 0$.
Thus, the assumptions on the marked Cox process $\X_M$ in Lemma
\ref{lemcondcox} can be reduced to $\E\|\Xi_0\| ^{d+\varepsilon} <
\infty$, respectively, $\E\|\Xi_0\|^{2d+\varepsilon} < \infty$ for arbitrarily
small $\varepsilon> 0$. The fact that $\X_M$ is $m$-dependent if
$\|\Xi_0\|$ is bounded allows us to apply an approximation technique
with truncated grains as in \cite{Hei99}, pages 299--302, showing that
the conditions with $\varepsilon= 0$ suffice. There exist substantial
examples of $\beta $-mixing PPs (e.g., obtained by dependent thinning or
clustering) which are far from being $m$-dependent. An example is
formed by the vertices of Poisson--Voronoi cells yielding exponentially
decaying $\beta$-mixing rates, see \cite{Hei94} for details.
\end{Remark*}

\section{Simulation study}\label{secstaapp}

Our aim was to find out whether the goodness-of-fit test for the Palm
mark distribution suggested by (\ref{teststat}) is suitable for the
detection of anisotropy
in Boolean models using directionally marked Cox processes on their
boundary as defined in Section \ref{seccoxBoolean}. This approach has
been applied to quality control of tomographic reconstruction
algorithms, see \cite{LueckMRS12}. Such algorithms typically introduce
elongation artifacts of objects when the input data suffers from a
missing wedge of projection angles as typical for electron tomography,
see \cite{midgley2003}. The accuracy of data varies locally with the
geometry of the specimen and may be reduced by use of appropriate
reconstruction algorithms, see \cite{LueckMRS12}. Our study is based
on simulated 2D Boolean models formed by discs with gamma distributed
radii (scale and shape parameter $4.5$ and $9$). These can be viewed as
2D slices of a 3D tomographic reconstruction of a complex foam-like
material. Note that in the parallel beam geometry of electron
tomography 3D volumes are stacks of 2D reconstructions generated from
1D projection data, which motivates this model choice in view of
the application in \cite{LueckMRS12}. Anisotropy artifacts were
simulated by transformation of the discs into ellipsoids with axes
parallel to the coordinate system. The major axis lengths were taken as
multiples of the minor axis lengths for factors $c_e\in\{1.135,1.325 \}
$. These values are typical elongation factors of standard
reconstruction algorithms for missing wedges of $30^\circ$ and
$60^\circ$, respectively, see \cite{LueckMRS12}.
The intensity of the Poisson PP ${\mathbf Y}$ of germs was chosen as
$1.5\cdot10^{-4}$ and the intensity of the Poisson PP of boundary
points as $0.1$.

Our asymptotic $\chi^2$-goodness-of-fit test is based on the test
statistic $T_k$ defined in (\ref{teststat}).
If $(P_M^{\mathbf o})_{0}$ denotes a hypothetical Palm mark
distribution, the hypothesis $H_0\dvt P_M^{\mathbf o}=(P_M^{\mathbf
o})_{0}$ is
rejected, if
$T_k>\chi^2_{\ell,1-\alpha}$, where $\alpha$ is the level of
significance, and $\chi^2_{\ell,1-\alpha}$ denotes the $(1-\alpha
)$-quantile of the $\chi^2_\ell$-distribution.
The bins
$C_1,\ldots, C_\ell\in\mathcal{B}({\mathbb S}_+^1)$ for the $\chi
^2$-goodness-of-fit test were chosen as
\[
C_i= \biggl\{(\cos\theta, \sin\theta)^T\dvt  \theta\in
\biggl[(i-1)\frac{\uppi }{\ell+1}, i\frac{\uppi }{\ell+1} \biggr) \biggr\},\qquad i=1,\ldots,\ell.
\]
We will discuss the case $\ell=8$, where
the bins had a width of $20^\circ$.
If $(\widehat{\covM})_k$ in (\ref{teststat}) is chosen as the
$L^2$-consistent estimator $(\widehat{\sigma}_{ij}^{(3)})_k$, the
test\vspace*{-1pt} will be referred to as
``\emph{test for the typical mark distribution}'' (TMD).
The construction of $(\widehat{\sigma}_{ij}^{(3)})_k$ involves the
sequence of bandwidths $\{b_k\}
$ chosen as
%
\begin{equation}\label{formbkkappa}
b_k=c |W_k|^{-{3}/({4d})} \qquad\mbox{for some constant
} c>0.
\end{equation}
The constant $c$ is crucial for test performance, as discussed below.
The asymptotic behavior of the tests was studied by considering squared
observation windows
corresponding to an expected number of $300$, $600,\ldots,3000$ points.
Due to the corresponding side lengths of the observation windows, (\ref
{formbkkappa}) entailed Condition \hyperref[condWB]{$(wb)$} and hence $(\widehat
{\sigma}_{ij}^{(3)})_k$ was
$L^2$-consistent.\vspace*{1pt}

The choice of the bandwidths $\{b_k\}$ can be avoided if $\covM$ is
not estimated
from the data to be tested but incorporated into $H_0$. This means, we
specify an MPP as null model, such that $\covM_0$ is either
theoretically known or otherwise can be approximated by Monte Carlo simulation.
By means of the combined null hypothesis $H_0\dvt  P_M^{\mathbf
o}=(P_M^{\mathbf o}
)_{0}$ and ${\covM}={\covM}_0$,
the test exploits
not only information on the distribution of the typical mark but
additionally considers asymptotic effects of
spatial dependence. The test can thus be used to investigate if a given
point pattern differs from the MPP null model w.r.t. the Palm mark
distribution. We will therefore refer to it as ``\emph{test for
mark-oriented goodness of model fit}''
(MGM). By the strong law of large numbers and the asymptotic
unbiasedness of $(\widehat{\sigma}_{ij}^{(2)})_k$,
a strongly consistent Monte Carlo estimator for $\covM_0$ in an MPP
model $\X_M$ is given by
\[
\widehat{\covM}_{k,n}=\frac{1}{n}\sum
_{\nu=1}^n \bigl(\widehat{\sigma
}_{ij}^{(2)}\bigr)_k\bigl(
\X_M^{(\nu)}\bigr),
\]
where $\X_M^{(1)},\ldots, \X_M^{(n)}$ are independent realizations
of $\X_M$.
Thus, for large $k$ and $n$ the test statistic $T_{k,n}={\mathbf Y}_k^\top
\widehat{\covM}_{k,n}^{-1}{\mathbf Y}_k $ has an approximate $\chi
^2_\ell$ distribution.
The estimator $\widehat{\covM}_{k,n}$ can also be used to construct a
test for the typical mark distribution if independent replications of a
point patterns are to be tested. In that case $\X_M^{(1)},\ldots, \X
_M^{(n)}$ are the replications. Note that for replicated point
patterns, $H_0$ does not incorporate an assumption on $\covM$ and
hence the corresponding test differs from the MGM test. The
edge-corrected unbiased estimator $(\widehat{\sigma}_{ij}^{(1)})_k$
was not used for the
Monte Carlo estimates in our simulation study, since $(\widehat{\sigma
}_{ij}^{(2)})_k$ can be
computed more efficiently.

All simulation results are based on $1000$ model realizations per
scenario. Type II errors were computed for
Boolean models with elongated grains, which means that the mark
distribution was not uniform on ${\mathbb S}_+^1$, whereas $H_0\dvt
P_M^{\mathbf o}
=U({\mathbb S}_+^1)$ hypothesized a uniform Palm mark distribution on
${\mathbb S}_+^1$.

\begin{figure}

\includegraphics{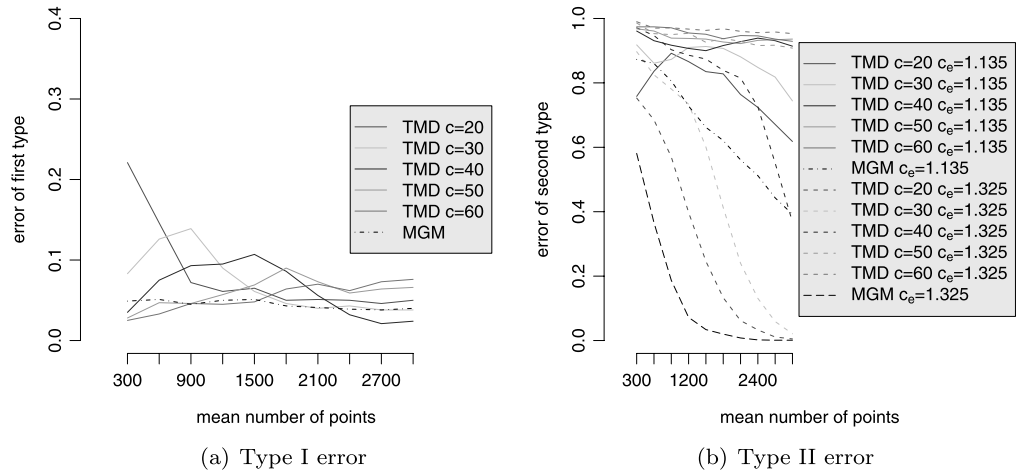}

\caption{Empirical errors of types I and II for the TMD and the MGM
test plotted against the mean number of points in the observation
window ($\alpha=0.5$). The constant $c$ is a bandwidth parameter for
covariance estimation in the TMD test, whereas $c_e$ denotes the
elongation factor of ellipses forming the Boolean model used as input
data for the analysis of type II errors.}
\label{figtestResults}
\end{figure}

The performance of the MGM test is visualized in Figure \ref
{figtestResults}. Empirical type I errors of the MGM test were close
to the theoretical $5\%$ level of significance, at which all tests were
conducted.
Experiments with the TMD test revealed that the choice of the bandwidth
parameter $c$ in (\ref{formbkkappa}) is critical for test
performance (Figure \ref{figtestResults}). Whereas large values of
$c$ result in a correct level of type I errors, they decrease the power
of the test. On the other hand, small values for $c$ lead to superior
power but at least for small observation windows with a limited number
of points increase type I errors (Figure \ref{figtestResults}).

The relatively high errors of second type for the small elongation
factor of $c_e=1.135$ are to be expected, since the investigated
structures are only slightly anisotropic. Nevertheless, for an expected
number of $3000$ points the MGM and TMD tests achieve a power of
$\sim60 \%$ and $40\%$, respectively, for $c_e=1.135$ and reject the
null hypothesis with probabilty $1$ for $c_e=1.325$. In summary, our
simulation results indicate that the MGM test outperforms the TMD test
especially with respect to power. This result is plausible since the
additional information incorporated into $H_0$ by specification of a
model covariance matrix can be expected to result in a more specific
test.

\section*{Acknowledgements}

We are grateful to the anonymous referees for their valuable
suggestions to improve the manuscript.




\printhistory

\end{document}